\documentclass[12pt,a4paper]{article}
\usepackage{amsmath,amsthm,amsfonts,graphicx,psfrag,amssymb,dsfont,mathrsfs}
\usepackage{minitoc,color}

\usepackage[colorlinks]{hyperref}
\usepackage[alphabetic,initials]{amsrefs}
\usepackage{bbm}
\usepackage{dsfont}
\usepackage[ngerman,english]{babel}
\usepackage{stackengine}
\usepackage{tikz}
\usepackage{enumitem}
\usepackage[ngerman,english]{babel}

\usetikzlibrary{matrix,arrows}

\newcommand{\A}{\mathbf{A}}

\newcommand{\ind}{\operatorname{ind}}

\theoremstyle{plain}
\newtheorem{thm}{Theorem}[section]

\newtheorem{prop}[thm]{Proposition}
\newtheorem{cor}[thm]{Corollary}
\newtheorem{lemma}[thm]{Lemma}

\theoremstyle{definition}
\newtheorem{example}[thm]{Example}
\newtheorem*{notation}{Notation}
\newtheorem{defn}[thm]{Definition}

\newtheorem{disclaimer}[thm]{Disclaimer}
\newtheorem{remark}[thm]{Remark}

\newtheorem{assumptions}[thm]{Assumptions}

\usepackage{scalerel}

\newcommand\reallywidehat[1]{\arraycolsep=0pt\relax%
\begin{array}{c}
\stretchto{
  \scaleto{
    \scalerel*[\widthof{\ensuremath{#1}}]{\kern-.5pt\bigwedge\kern-.5pt}
    {\rule[-\textheight/2]{1ex}{\textheight}} 
  }{\textheight} %
}{0.8ex}\\           
#1\\                 
\rule{-1ex}{0ex}
\end{array}
}

\date{}
\title{Algebraic and Giroux torsion in higher-dimensional contact manifolds}
\author{Agustin Moreno}

\begin{document}

\maketitle
\thispagestyle{empty}

\begin{abstract} We construct examples in any odd dimension of contact manifolds with finite and non-zero algebraic torsion (in the sense of \cite{LW}), which are therefore tight and do not admit strong symplectic fillings. We prove that Giroux torsion implies algebraic $1$-torsion in any odd dimension, which proves a conjecture in \cite{MNW}. These results are part of the author's PhD thesis \cite{Mo2}.
\end{abstract} 

\addtocontents{toc}{\protect\sloppy}
\tableofcontents

\section{Introduction}

In this paper, and its followup \cite{Mo}, we address the general problem of constructing ``interesting'' examples of higher-dimensional contact manifolds,  and developing techniques in order to compute SFT-type holomorphic curve invariants. 

We will construct examples of contact manifolds in every odd dimension, presenting a geometric structure which is a higher-dimensional version of that of a \emph{spinal open book decomposition} or SOBD, as defined in \cite{LVHMW} in dimension $3$. The type of SOBD present in our examples, which one could call \emph{partially planar}, mimics the notion of \emph{planar $m$-torsion} domains as defined in \cite{Wen2}. Indeed, it consists of two surface fibrations over a higher-dimensional contact base, one of them having genus zero fibers, glued together along a contact fibration over a Liouville domain. This geometric structure can be ``detected'' algebraically by algebraic torsion, a holomorphic-curve contact invariant. For suitable data, the surface fibers become holomorphic, and are leaves of a finite energy foliation of the symplectization $\mathbb{R}\times M$. The isolated ones may be counted in a suitable way, and the result is an invariant which ``recovers'' the number $m$. This is the idea inspiring algebraic $m$-torsion.

We exhibit a detailed construction of an isotopy class of contact forms, which is ``supported'' by the SOBD, so that one may view these contact forms as ``Giroux'' forms. We will estimate the algebraic torsion of these examples, which we show is finite, and, in certain cases, non-zero. In those cases, the contact manifolds are tight and admit no strong symplectic fillings. 

We will also relate algebraic torsion with a geometric condition, \emph{Giroux torsion}. While this is a classical notion in dimension 3, the higher-dimensional version was introduced in \cite{MNW}. We will show that the geometric presence of certain torsion domains inside a contact manifold can be detected algebraically by SFT. More concretely, Giroux torsion implies algebraic $1$-torsion, in any odd dimension. This proves a conjecture in \cite{MNW}. 

The proof of this result is carried out by interpreting the Giroux torsion domains as being supported by a suitable SOBD, which we call a \emph{Giroux SOBD}, for which we give a notion of a ``Giroux form''. The result follows by adapting our computations for the above partially planar model contact manifolds. 

In order to carry out our computations, we need a very detailed understanding of holomorphic curves in the symplectization of our model contact manifolds, and the SOBD structures, together with the associated finite energy foliations, are crucial towards this end. The key technical inputs are: transversality of the genus zero curves in the foliation, and a uniqueness result for holomorphic curves (Theorem \ref{uniqueness}). Proving transversality is needed so that indeed one has a space of isolated curves to count, whereas uniqueness is necessary to know precisely what to count. For transversality, a standard technique in dimension three is the automatic transversality criterion of \cite{Wen1}, which consists in checking a fairly straightforward numerical inequality involving topological data associated to a given curve. For uniqueness, one can sometimes resort to Siefring's intersection theory for punctured holomorphic curves in dimension four \cite{Sie11}. In higher dimensions, things become cumbersome. To prove transversality, we resorted to a ``hands-on'' analytical approach of computing precisely the kernel of the linearization of the Cauchy-Riemann operator, and check that its dimension coincides with its Fredholm index, from which transversality follows. For uniqueness, we resorted to a combination of energy estimates, holomorphic cascades and geometric arguments.

\subparagraph*{On the invariant.} The invariant we will use, algebraic torsion, was defined in \cite{LW}, and is a contact invariant taking values in $\mathbb{Z}^{\geq 0}\cup \{\infty\}$. It was introduced, using the machinery of \emph{Symplectic Field Theory}, as a quantitative way of measuring non-fillability, giving rise to a ``hierarchy of fillability obstructions'', cf.\ \cite{Wen2}. At least morally, $0$-torsion should correspond to overtwistedness, whereas $1$-torsion is implied by Giroux torsion (the converse is not true). Having $0$-torsion is actually equivalent to being \emph{algebraically overtwisted}, which means that the contact homology, or equivalently its SFT, vanishes (Proposition 2.9 in \cite{LW}). This is well-known to be implied by overtwistedness, but the converse is still wide open.

The key fact about this invariant is that it behaves well under exact symplectic cobordisms, which implies that the concave end inherits any order of algebraic torsion that the convex end has. Thus, algebraic torsion may be also thought of as an obstruction to the existence of exact symplectic cobordisms. In particular, it serves as an obstruction to symplectic fillability. Moreover, there are connections to dynamics: any contact manifold with finite torsion satisfies the Weinstein conjecture (i.e.\ there exist closed Reeb orbits for every contact form).  

One should mention that there are other notions of algebraic torsion in the literature which do not use SFT, but which are only 3-dimensional (see \cite{KMvhMW} for the version using Heegard Floer homology, or the appendix in \cite{LW} by Hutchings, using ECH).

\subparagraph*{Statement of results.}

For the SFT setup, we follow \cite{LW}, where we refer the reader for more details. We will take the SFT of a contact manifold $(M,\xi)$ (with coefficients) to be the homology $H^{SFT}_*(M, \xi; \mathcal{R})$ of a $\mathbb{Z}_2$-graded unital $BV_\infty$-algebra $(\mathcal{A}[[\hbar]], \mathbf{D}_{SFT})$ over the group ring $R_{\mathcal{R}}:=\mathbb{R}[H_2(M;\mathbb{R})/\mathcal{R}]$, for some linear subspace $\mathcal{R}\subseteq H_2(M;\mathbb{R})$. Here, $\mathcal{A}=\mathcal{A}(\lambda)$ has generators $q_\gamma$ for each good closed Reeb orbit $\gamma$ with respect to some nondegenerate contact form $\lambda$ for $\xi$, $\hbar$ is an even variable, and the operator $$\mathbf{D}_{SFT}: \mathcal{A}[[\hbar]] \rightarrow\mathcal{A}[[\hbar]]$$ is defined by counting rigid solutions to a suitable abstract perturbation of a $J$-holomorphic curve equation in the symplectization of $(M, \xi)$. It satisfies

\begin{itemize}
 \item $\mathbf{D}_{SFT}$ is odd and squares to zero,
 \item $\mathbf{D}_{SFT}(1) = 0$, and
 \item $\mathbf{D}_{SFT} = \sum_{k \geq 1} D_k\hbar^{k-1},$
\end{itemize}

where $D_k : \mathcal{A} \rightarrow \mathcal{A}$ is a differential operator of order $\leq k$, given by
$$
D_k = \sum_{\substack{\Gamma^+,\Gamma^-, g,d\\ |\Gamma^+|+g=k}}\frac{n_g(\Gamma^+,\Gamma^-,d)}{C(\Gamma^-,\Gamma^+)}q_{\gamma_1^-}\dots q_{\gamma_{s^-}^-}z^d\frac{\partial}{\partial q_{\gamma_1^+}}\dots \frac{\partial}{\partial q_{\gamma_{s^+}^+}}
$$
The sum ranges over all non-negative integers $g \geq 0$, homology classes $d \in H_2(M; \mathbb{R})/\mathcal{R}$ and ordered (possibly empty) collections of good closed Reeb orbits $\Gamma^\pm = (\gamma_1^{\pm},\dots,\gamma_{s^\pm}^\pm)$ such that $s^+ + g = k$. After a choice of spanning surfaces as in \cite{EGH} (p. 566, see also p. 651), the projection to $M$ of each finite
energy holomorphic curve $u$ can be capped off to a 2-cycle in $M$, and so it gives rise to a homology class $[u]\in H_2(M)$, which we project to define $\overline{[u]} \in H_2(M; \mathbb{R})/\mathcal{R}$. The number $n_g(\Gamma^+,\Gamma^-,d) \in \mathbb{Q}$ denotes the count of (suitably perturbed) holomorphic curves of genus $g$ with positive asymptotics $\Gamma^+$ and negative asymptotics $\Gamma^-$ in the homology class $d$, including asymptotic markers as explained in \cite{EGH}, or \cite{Wen3}, and including rational weights arising from automorphisms. $C(\Gamma^-, \Gamma^+) \in \mathbb{N}$ is a combinatorial factor defined as $C(\Gamma^-, \Gamma^+) = s^-!s^+!\kappa_{\gamma_1^-}\dots\kappa_{\gamma_{s^-}^-}$, where $\kappa_\gamma$ denotes the covering multiplicity of the Reeb orbit $\gamma$.

The most important special cases for our choice of linear subspace $\mathcal{R}$ are $\mathcal{R} = H_2(M; \mathbb{R})$ and $\mathcal{R} = \{0\}$, called the \emph{untwisted} and \emph{fully twisted} cases respectively, and $\mathcal{R} = \ker \Omega$ with $\Omega$ a closed 2-form on $M$. We shall abbreviate the latter case as $H^{SFT}_*(M, \xi;\Omega) := H^{SFT}_*(M, \xi; \ker \Omega)$, and the untwisted case simply by $H^{SFT}_*(M, \xi):=H^{SFT}_*(M, \xi;H_2(M;\mathbb{R}))$.

\begin{defn} Let $(M, \xi)$ be a closed manifold of dimension $2n+1$ with a positive, co-oriented contact structure. For any integer $k \geq 0$, we say that $(M, \xi)$ has $\Omega$-twisted algebraic torsion of order $k$ (or $\Omega$-twisted $k$-torsion) if $[\hbar^k] = 0$ in $H^{SFT}_*(M, \xi;\Omega)$. If this is true for all $\Omega$, or equivalently, if $[\hbar^k] = 0$ in $H^{SFT}_*(M, \xi;\{0\})$, then we say that $(M, \xi)$ has fully twisted algebraic $k$-torsion. 
\end{defn}

We will refer to \emph{untwisted} $k$-torsion to the case $\Omega=0$, in which case $R_\mathcal{R}=\mathbb{R}$ and we do not keep track of homology classes. Whenever we refer to torsion without mention to coefficients we will mean the untwisted version. We will say that, if a contact manifold has algebraic $0$-torsion for every choice of coefficient ring, then it is \emph{algebraically overtwisted}, which is equivalent to the vanishing of the SFT, or its contact homology. By definition, $k$-torsion implies $(k+1)$-torsion, so we may define its \emph{algebraic torsion} to be
$$
AT(M,\xi;\mathcal{R}):=\min\{k\geq 0:[\hbar^k]=0\}\in \mathbb{Z}^{\geq 0}\cup\{\infty\},
$$
where we set $\min \emptyset=\infty$. We denote it by $AT(M,\xi)$, in the untwisted case.

This construction is well-behaved under symplectic cobordisms: 
Any exact symplectic cobordism $(X, \omega=d\alpha)$ with positive end $(M_+, \xi_+)$ and negative end $(M_-, \xi_-)$ gives rise to a natural $\mathbb{R}[[\hbar]]$-module morphism on the untwisted SFT,
$$\Phi_X: H^{SFT}_*(M_+, \xi_+) \rightarrow H^{SFT}_*(M_-, \xi_-),$$ a \emph{cobordism map}. This implies that if $(M_+, \xi_+)$ has $k$-torsion, then so does $(M_-, \xi_-)$. There is also a version with coefficients for the case of non-exact cobordisms and fillings \cite[Prop.\ 2.4]{LW}.

Examples of 3-dimensional contact manifolds with any given order of torsion $k-1$, but not $k-2$, were constructed in \cite{LW}. The underlying manifold is the product manifold $M_g:=S^1\times \Sigma$, for $\Sigma$ a surface of genus $g$ which is divided into two pieces $\Sigma_+$ and $\Sigma_-$ along some dividing set of simple closed curves $\Gamma$ of cardinality $k$, where the latter has genus $0$, and the former has genus $g-k+1$. The contact structure $\xi_k$ is $S^1$-invariant and may be obtained, for instance, by a construction originally due to Lutz (see \cite{Lutz}). Its isotopy class is characterized by the fact that every section $\{pt\}\times\Sigma$ is a convex surface with dividing set $\Gamma$. The behaviour of algebraic torsion under cobordisms then implies that there is no exact symplectic cobordisms having $(M_g,\xi_k)$ and $(M_{g^\prime},\xi_{k^\prime})$ as convex and concave ends, respectively, if $k<k^\prime$.

The existence of the analogue higher dimensional contact manifolds was conjectured in \cite{LW}. We will consider a modified version of their examples. The modification we do here consists in taking the $S^1$-factor and replacing it by a \emph{closed} $(2n-1)$-manifold $Y$, having the special property that $Y\times I$ admits the structure of a Liouville domain (here, $I$ denotes the interval $[-1,1]$). This means that it comes with an exact symplectic form $d\alpha$, and has disconnected contact-type boundary $\partial(Y \times I,d\alpha)=(Y_-,\xi_-=\ker\alpha_-)\bigsqcup (Y_+,\xi_+=\ker\alpha_+)$, where $Y_\pm$ coincide with $Y$ as manifolds, but $(Y_\pm,\xi_\pm)$ are \emph{not contactomorphic} to each other. In fact, $Y_\pm$ have different orientations, and so they might not even be \emph{homeomorphic} to each other (not every manifold admits an orientation-reversing homeomorphism). A Liouville domain of the form $(Y\times I,d\alpha)$ is what we will call a \emph{cylindrical Liouville semi-filling} (or simply a cylindrical semi-filling). Their existence in every odd dimension was established in \cite{MNW}. We immediately see that this generalizes the previous 3-dimensional example, since $S^1$ admits the \emph{Liouville pair} $\alpha_\pm=\pm d\theta$, which means that the 1-form $e^{-s}\alpha_- +e^s\alpha_+$ is Liouville in $S^1 \times \mathbb{R}$. We prove that the manifold $Y\times \Sigma$ indeed achieves $(k-1)$-torsion (Theorem \ref{thm5d}), for a suitable contact structure which we now describe. 

First, for once and for good, we will fix the following notation:

\begin{notation}\label{not} Throughout this paper, the symbol $I$ will be reserved for the interval $[-1,1]$.
\end{notation}

\begin{figure}[t]\centering \includegraphics[width=0.70\linewidth]{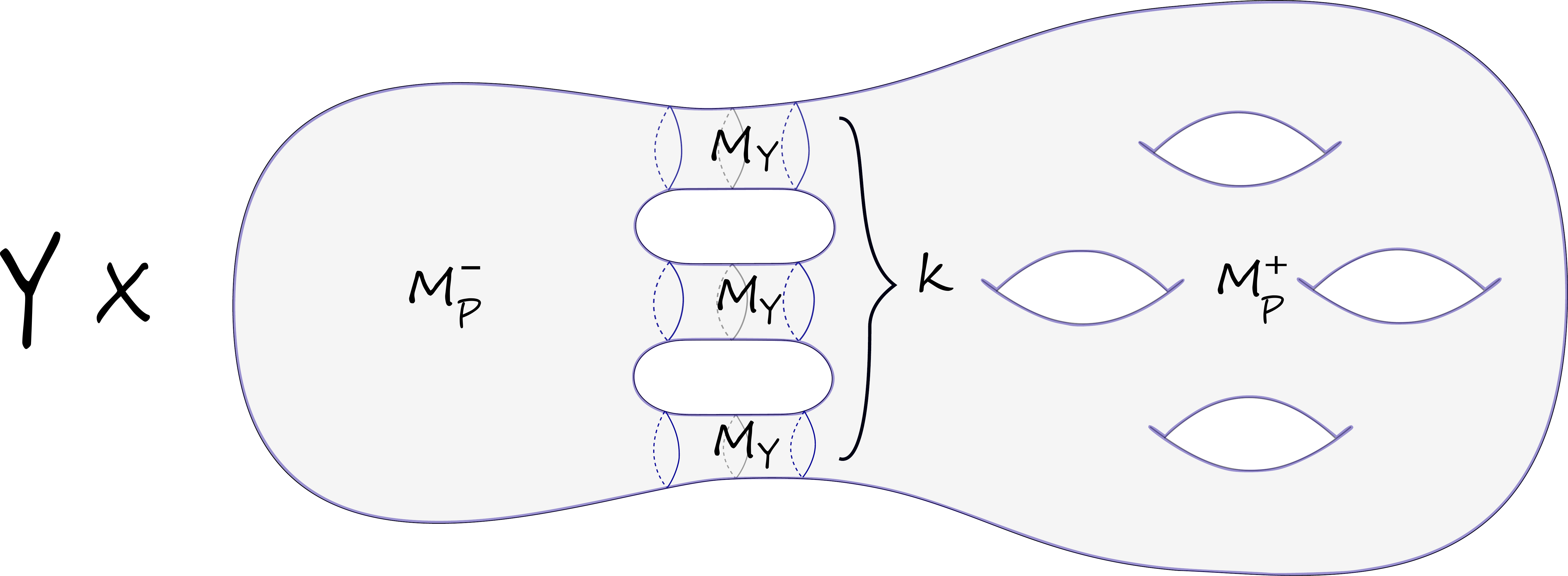}
\caption{\label{themanifold} The SOBD structure in $M$.}
\end{figure}

We can adapt the construction of the contact structures in \cite{LW} to our models. The starting idea is to decompose the manifold $	M=M_g=Y\times \Sigma$ into three pieces
$$
M_g=M_Y \bigcup M_P^\pm,
$$
where $M_Y = \bigsqcup^k Y \times I \times S^1$, and $M_P^\pm = Y \times \Sigma_\pm$ (see Figure \ref{themanifold}). 
We have natural fibrations $$\pi_Y: M_Y \rightarrow Y \times I$$$$ \pi_P^\pm: M_P^{\pm} \rightarrow Y_\pm,$$ with fibers $S^1$ and $\Sigma_\pm$, respectively, and they are compatible in the sense that
$$
\partial((\pi_P^\pm)^{-1}(pt))=\bigsqcup^k \pi_Y^{-1}(pt)
$$

While $\pi_Y$ has a Liouville domain as base, and a contact manifold as fiber, the situation is reversed for $\pi_P^\pm$, which has contact base, and Liouville fibers. This is a prototypical example of a \emph{spinal open book decomposition}, or SOBD. While we will not give a general definition of such a notion, we refer the reader to \cite{Mo2} for a tentative one. 

Using this decomposition, we can construct a contact structure $\xi_k$ which is a small perturbation of the stable Hamiltonian structure $\xi_\pm \oplus T\Sigma_\pm$ along $M_P^\pm$, and is a \emph{contactization} for the Liouville domain $(Y\times I,\epsilon d\alpha)$ along $M_Y$, for some small $\epsilon>0$. This means that it coincides with $\ker(\epsilon\alpha + d\theta)$, where $\theta$ is the $S^1$-coordinate. We will do this in detail in Section \ref{modelAsec}.

\begin{remark}\label{dualSOBD} Let us remark that, since the fibrations above are trivial, one can always reverse their roles. More precisely, we could consider instead the ``dual'' SOBD:
$$\pi^*_Y: M_Y \rightarrow S^1$$$$ (\pi_P^\pm)^*: M_P^{\pm} \rightarrow \Sigma_\pm$$
For these fibrations, we may also construct a contact form which is ``supported'' by the SOBD. The resulting contact structure is isotopic to $\xi_k$, which is what we expect from the point of view of a ``Giroux correspondence'' (in this more general setting). This is actually used for the results in \cite{Mo}. 
\end{remark}

For the contact manifolds $(M_g,\xi_k)$, we can estimate their algebraic torsion. First, recall that a contact structure is hypertight if it admits a contact form without contractible Reeb orbits (which we call a hypertight contact form). In particular, there are no holomorphic disks in their symplectization, which implies that there is no 0-torsion. By a well-known theorem by Hofer and its generalization to higher dimensions by Albers--Hofer (in combination with \cite{BEM}), hypertight contact manifolds are tight.  

\begin{thm}\label{thm5d}
For any $k\geq 1$, and $g\geq k$, the $(2n+1)$-dimensional contact manifolds $(M_g=Y\times \Sigma,\xi_k)$ satisfy $AT(M_g,\xi_k)\leq k-1$. Moreover, if $(Y,\alpha_\pm)$ are hypertight, and $k\geq 2$, the corresponding contact manifold $(M_g,\xi_k)$ is also hypertight. In particular, $AT(M_g,\xi_k)>0$, and it is tight.  
\end{thm}

In fact, the examples of Theorem \ref{thm5d} admit $\Omega$-twisted $k-1$-torsion, for $\Omega$ defining a cohomology class in $\mathcal{O}:=\mbox{Ann}(\bigoplus_k H_1(Y;\mathbb{R})\otimes H_1(S^1; \mathbb{R}))$, the annihilator of $\bigoplus_k H_1(Y;\mathbb{R})\otimes H_1(S^1; \mathbb{R}) \subseteq H_2(M_g;\mathbb{R})$. Here, we take the homology of the subregion $\bigsqcup^k Y\times \{0\}\times S^1$, lying along the region $M_Y$ where $\Sigma_\pm$ glue together. Using \cite[Prop. 2.4]{LW}, we obtain:

\begin{cor}\label{thm5d2}
The examples of Theorem \ref{thm5d} do not admit weak fillings $(W,\omega)$ for which $[\omega\vert_{M_g}]$ is rational and lies in $\mathcal{O}$. In particular, they are not strongly fillable. 
\end{cor} 

\begin{remark} \label{mit}
\begin{itemize}[wide, labelwidth=!, labelindent=0pt] $\;$
\item By a result of Mitsumatsu in \cite{Mit}, any 3-manifold $Y$ which admits a smooth Anosov flow preserving a smooth volume form satisfies that $Y \times I$ can be enriched with a cylindrical Liouville semi-filling structure. Therefore any of these $3$-manifolds can be used in the construction of $5$-dimensional contact models with $AT\leq k-1$, for any $k\geq 1$.

\item The examples of Liouville cylindrical semi-fillings of \cite{MNW} satisfy the hypertightness condition. Then we have a doubly-infinite family of contact manifolds with $0< AT(M_g,\xi_k)\leq k-1$, in any dimension. These are then an instance of higher-dimensional tight but not strongly fillable contact manifolds, since they have non-zero and finite algebraic torsion. For $k=2$, this precisely computes the algebraic torsion.
\end{itemize}
\end{remark}

The authors of \cite{MNW} define a generalized higher-dimensional version of the notion of Giroux torsion. This notion is defined as follows: consider $(Y,\alpha_+,\alpha_-)$ a \emph{Liouville pair} on a closed manifold $Y^{2n-1}$, which means that the $1$-form $\beta=\frac{1}{2}(e^s\alpha_+ + e^{-s}\alpha_-)$ is Liouville in $\mathbb{R}\times Y$. Consider also the \emph{Giroux $2\pi$-torsion domain} modeled on $(Y,\alpha_+,\alpha_-)$ given by the contact manifold $(GT,\xi_{GT}):=(Y\times [0,2\pi]\times S^1,\ker\lambda_{GT})$, where 
\begin{equation}
\begin{split}
\lambda_{GT}&=\frac{1+\cos(r)}{2}\alpha_++\frac{1-\cos(r)}{2}\alpha_-+\sin (r)d\theta\\
\end{split}
\end{equation} and the coordinates are $(r,\theta)\in [0,2\pi] \times S^1$. Say that a contact manifold $(M^{2n+1},\xi)$ has \emph{Giroux torsion} whenever it admits a contact embedding of $(GT,\xi_{GT})$. In this situation, denote by $\mathcal{O}(GT)\subseteq H^2(M;\mathbb{R})$ the annihilator of $\mathcal{R}_{GT}:=H_1(Y;\mathbb{R})\otimes H_1(S^1;\mathbb{R})$, viewed as a subspace of $H_2(M;\mathbb{R})$. The following was conjectured in \cite{MNW}:

\begin{thm}\label{AlgvsGiroux}
If a contact manifold $(M^{2n+1},\xi)$ has Giroux torsion, then it has $\Omega$-twisted algebraic 1-torsion, for every $[\Omega] \in \mathcal{O}(GT)$, where $GT$ is a Giroux $2\pi$-torsion domain embedded in $M$. 
\end{thm}

The proof uses the same techniques as Theorem \ref{thm5d}, and the main idea is to interpret Giroux torsion domains in terms of a specially simple kind of SOBD, which we call \emph{Giroux SOBD}. 

\vspace{0.5cm}

A natural corollary is the following:

\begin{cor}\label{GTfill} If a contact manifold $(M^{2n+1},\xi)$ has Giroux torsion, then it does not admit weak fillings $(W,\omega)$ with $[\omega\vert_M]\in \mathcal{O}(GT)$ and rational, where $GT$ is a Giroux $2\pi$-torsion domain embedded in $M$. In particular, it is not strongly fillable.
\end{cor}

This is essentially corollary 8.2 in \cite{MNW}, which was obtained with different methods. Observe that if $\mathcal{R}_{GT}=0$ then $(M,\xi)$ does not admit weak fillings at all. This is in fact the condition used in \cite{MNW} to obstruct weak fillability.

\subparagraph*{Further work: a synopsis.} We now state a series of results, to be proven in the followup paper \cite{Mo} (see also \cite{Mo2}). In the following, we use the fact that the unit cotangent bundle of a hyperbolic surface fits into a cylindrical semi-filling \cite{McD}.

\begin{thm}\label{nocob1} Let $(M_0^{5},\xi_0)$ be a $5$-dimensional contact manifold with Giroux torsion, and let $Y$ be the unit cotangent bundle of a hyperbolic surface. If $(M=\Sigma \times Y,\xi_k)$ is the corresponding $5$-dimensional contact manifold of Theorem \ref{thm5d} with $k\geq 3$, then there is no exact symplectic cobordism having $(M_0,\xi_0)$ as the convex end, and $(M,\xi_k)$ as the concave end.
\end{thm}

In particular, we obtain

\begin{cor}\label{corGT}
If $Y$ is the unit cotangent bundle of a hyperbolic surface, and $(M=\Sigma \times Y,\xi_k)$ is the corresponding $5$-dimensional contact manifold of Theorem \ref{thm5d} with $k\geq 3$, then $(M,\xi_k)$ does not have Giroux torsion.
\end{cor}

Moreover, we have reasons, coming from string topology \cite{CL09}, to believe that the examples of Corollary \ref{corGT} have untwisted algebraic 1-torsion (for \emph{any} $k\geq 1$).

Putting Theorem \ref{thm5d} (and Remark \ref{mit}), together with Corollaries \ref{thm5d2} and \ref{corGT}, we obtain the following:

\begin{cor}\label{corG}
There exist infinitely many non-diffeomorphic $5$-dimensional contact manifolds $(M,\xi)$ which are tight, not strongly fillable, and which do not have Giroux torsion.
\end{cor}

To our knowledge, there are no other known examples of higher-dimensional contact manifolds as in Corollary \ref{corG}. Also, we expect the above examples to have algebraic 1-torsion. 

\vspace{0.5cm}

One can twist the contact structure of Theorem \ref{thm5d} close to the dividing set, by performing the \emph{$l$-fold Lutz--Mori twist} along a hypersurface $H$ lying in $\partial(\bigsqcup^k Y \times I \times S^1)$. This notion was defined in \cite{MNW}, and builds on ideas by Mori in dimension 5 \cite{Mori09}. The resulting contact structures are, in general, all homotopic as almost contact structures, but in our case they are distinguishable by a suitable version of cylindrical contact homology. By construction, all of these have Giroux torsion, so by Theorem \ref{AlgvsGiroux} they have $\Omega$-twisted $1$-torsion, for $[\Omega] \in \mathcal{O}=\mathcal{O}(GT)$.

As a corollary of Theorem \ref{nocob1}, we get:

\begin{cor}\label{nocob}
Let $Y$ be the unit cotangent bundle of a hyperbolic surface, and let $(M=\Sigma \times Y,\xi)$ be the corresponding $5$-dimensional contact manifold of Theorem \ref{thm5d}, with $k\geq 3$. If $(M,\xi_l)$ denotes the contact manifold obtained by an $l$-fold Lutz--Mori twist of $(M,\xi)$, then there is no exact symplectic cobordism having $(M,\xi_l)$ as the convex end, and $(M,\xi)$ as the concave end (even though the underlying manifolds are diffeomorphic, and the contact structures are homotopic as almost contact structures).  
\end{cor}

The results from \cite{Mo} stated above make use of Richard Siefring's intersection theory for holomorphic curves and hypersurfaces, as outlined in an appendix in \cite{Mo2} written in coauthorship with Siefring, as a prequel of his upcoming work \cite{Sie2}, and to appear as an independent article \cite{MS}. Another technical input is the obstruction bundle technique as in Hutchings-Taubes \cite{HT1,HT2}. The SOBD is ``dualized'' in the sense of Remark \ref{dualSOBD}, and the finite energy foliation is replaced by a foliation by holomorphic hypersurfaces. Siefring's intersection theory then implies that holomorphic curves with suitable asymptotic behaviour lie in the leaves of the foliation. This, combined with symmetries in the setup and the obstruction bundle technique, allows us to obtain our results, as well as information on the SFT of our contact manifolds. 

\begin{disclaimer}\label{disc} Since the statements of our results make use of machinery from Symplectic Field Theory, they come with the standard disclaimer that they assume that its analytic foundations are in place. They depend on the abstract perturbation scheme promised by the polyfold theory of Hofer--Wysocki--Zehnder. We shall assume that it is possible to achieve transversality by introducing an arbitrarily small abstract perturbation to the Cauchy-Riemann equation, and that the analogue of the SFT compactness theorem still holds as the perturbation is turned off. In practice, this means that, in order to study curves for the perturbed data, we need to also study holomorphic building configurations for the unperturbed one. However, we have taken special care in that the approach taken not only provides results that will be fully rigorous after the polyfold machinery is complete, but also gives several direct results that are \emph{already} rigorous.
\end{disclaimer}

\subparagraph*{Acknowledgements.} 

First of all, my thanks go to my PhD supervisor, Chris Wendl, for introducing me to this project and for his support and patience throughout its duration. To Richard Siefring, for very helpful conversations and for co-authoring an appendix in \cite{Mo2}. To Janko Latschev and Kai Cieliebak, for going through the long process of reading \cite{Mo2}. To Patrick Massot, Sam Lisi, and Momchil Konstantinov, for helpful conversations/correspondence on different topics. 

This research, forming part of the author's PhD thesis, has been partly carried out in (and funded by) University College London (UCL) in the UK, and by the Berlin Mathematical School (BMS) in Germany.

\section*{Guide to the document}

The main construction is dealt with in Section \ref{modelAsec}. We show Fredholm regularity in Section \ref{Regularity}, and uniqueness (Theorem \ref{uniqueness}) in Section \ref{uniquenessMB}. Theorem \ref{thm5d} is proved in Section \ref{Torsion}.

The proof of Theorem \ref{AlgvsGiroux} is dealt with in Section \ref{GirouxTorsion1-tors}, which is basically a reformulation of the previous sections, with the key input being an adaptation of the uniqueness Theorem \ref{uniqueness}.
  
\newpage

\paragraph*{Basic notions}\label{basicn}

A contact form in a $(2n+1)$-dimensional manifold $M$ is a $1$-form $\alpha$ such that $\alpha \wedge d\alpha^{n}$ is a volume form, and the associated contact structure is $\xi=\ker \alpha$ (we will assume all our contact structures are co-oriented). The Reeb vector field associated to $\alpha$ is the unique vector field $R_\alpha$ on $M$ satisfying $$\alpha(R_\alpha)=1,\;i_{R_\alpha}d\alpha=0$$

A $T$-periodic Reeb orbit is $(\gamma,T)$ where $\gamma: \mathbb{R}\rightarrow M$ is such that $\dot{\gamma}(t)=TR_\alpha(\gamma(t))$, $\gamma(1)=\gamma(0)$. We will often just talk about a Reeb orbit $\gamma$ without mention to $T$, called its period, or action. If $\tau>0$ is the minimal number for which $\gamma(\tau)=\gamma(0)$, and $k \in \mathbb{Z}^+$ is such that $T=k\tau$, we say that the covering multiplicity of $(\gamma,T)$ is $k$. If $k=1$, then $\gamma$ is said to be simply covered (otherwise it is multiply covered). A periodic orbit $\gamma$ is said to be non-degenerate if the restriction of the time $T$ linearised Reeb flow $d\varphi^T$ to $\xi_{\gamma(0)}$ does not have $1$ as an eigenvalue. More generally, a Morse--Bott submanifold of $T$-periodic Reeb orbits is a closed submanifold $N \subseteq M$ invariant under $\varphi^T$ such that $\ker(d\varphi^T - \mathds{1})=TN$, and $\gamma$ is Morse--Bott whenever it lies in a Morse--Bott submanifold, and its minimal period agrees with the nearby orbits in the submanifold. The vector field $R_\alpha$ is non-degenerate/Morse--Bott if all of its closed orbits are non-degenerate/Morse--Bott.

A stable Hamiltonian structure (SHS) on $M$ is a pair $\mathcal{H}=(\Lambda,\Omega)$ consisting of a closed $2$-form $\Omega$ and a $1$-form $\Lambda$ such that 

$$\ker \Omega \subseteq \ker d\Lambda,\;\mbox{and}\;\Omega\vert_{\xi} \mbox{ is non-degenerate, where } \xi=\ker \Lambda$$ 

In particular, $(\alpha,d\alpha)$ is a SHS whenever $\alpha$ is a contact form. The Reeb vector field associated to $\mathcal{H}$ is the unique vector field on $M$ defined by $$\Lambda(R)=1,\; i_{R}\Omega=0$$ There are analogous notions of non-degeneracy/Morse--Bottness for SHS.

A symplectic form in a $2n$-dimensional manifold $W$ is a $2$-form $\omega$ which is closed and non-degenerate. A Liouville manifold (or an exact symplectic manifold) is a symplectic manifold with an exact symplectic form $\omega=d\lambda$, and the associated Liouville vector field $V$ is defined by the equation $i_{V}d\lambda=\lambda$. Any Liouville manifold is necessarily open. A boundary component $M$ of a Liouville manifold (endowed with the boundary orientation) is convex if the Liouville vector field is positively transverse to $M$, and is concave, if it is so negatively. An exact cobordism from a (co-oriented) contact manifold $(M_+,\xi_+)$ to $(M_-,\xi_-)$ is a compact Liouville manifold $(W,\omega=d\alpha)$ with boundary $\partial W=M_+\bigsqcup M_-$, where $M_+$ is convex, $M_-$ is concave, and $\ker \alpha\vert_{M_\pm}=\xi_\pm$. Therefore, the boundary orientation induced by $\omega$ agrees with the contact orientation on $M_+$, and differs on $M_-$. A Liouville filling (or a Liouville domain) of a --possibly disconnected-- contact manifold $(M,\xi)$ is a compact Liouville cobordism from $(M,\xi)$ to the empty set. A strong symplectic cobordism and a strong filling are defined in the same way, with the difference that $\omega$ is exact only in a neighbourhood of the boundary of $W$ (so that the Liouville vector field is defined in this neighbourhood, but not necessarily in its complement).

The symplectization of a contact manifold $(M,\xi=\ker \alpha)$ is the symplectic manifold $(\mathbb{R}\times M, \omega=d(e^a\alpha))$, where $a$ is the $\mathbb{R}$-coordinate. In particular, it is a non-compact Liouville manifold. Similarly, the symplectization of a stable Hamiltonian manifold $(M,\Lambda,\Omega)$ is the symplectic manifold $(\mathbb{R}\times M, \omega^\varphi)$, where $\omega^\varphi=d(\varphi(a)\Lambda)+\Omega$, and $\varphi$ is an element of the set $$\mathcal{P}=\{\varphi \in C^\infty(\mathbb{R},(-\epsilon,\epsilon)):\varphi^\prime>0\}$$ Here, $\epsilon>0$ is chosen small enough so that $\omega^\varphi$ is indeed symplectic. An $\mathcal{H}$-compatible (or simply cylindrical) almost complex structure on a symplectization $(W=\mathbb{R}\times M,\omega^\varphi)$ is $J \in \mbox{End}(TW)$ such that $$J \mbox{ is } \mathbb{R}\mbox{-invariant}, J^2=-\mathds{1},\;J(\partial_a)=R,\;J(\xi)=\xi,\;J\vert_{\xi} \mbox{ is }\Omega \mbox{-compatible }$$ The last condition means that $\Omega(\cdot,J\cdot)$ defines a $J$-invariant Riemannian metric on $\xi$. If $J$ is $\mathcal{H}$-compatible, then it is easy to check that it is $\omega^\varphi$-compatible, which means that $\omega^\varphi(\cdot,J\cdot)$ is a $J$-invariant Riemannian metric on $\mathbb{R}\times M$.

To any closed $T$-periodic Reeb orbit $(\gamma,T)$ one can associate an asymptotic operator $\mathbf{A}_\gamma$. To write it down, choose a symmetric connection $\nabla$ on $M$, and a $\mathcal{H}$-compatible almost complex structure $J$, and define
$$
\mathbf{A}_\gamma=\mathbf{A}_{\gamma,J}:W^{1,2}(\gamma^*\xi)\rightarrow L^2(\gamma^*\xi)
$$
$$
\mathbf{A}_\gamma\eta =-J(\nabla_t \eta -T\nabla_\eta R)
$$

Alternatively, one has the expression 
$$
\A_{\gamma}\eta(t)=-J\left.\frac{d}{ds}\right|_{s=0}d\varphi^{-Ts}\eta(t+s),
$$
for $\eta \in W^{1,2}(\gamma^*\xi)$, where $\varphi^s$ again is the time-$s$ Reeb flow.

Morally, this is the Hessian of a certain action functional on the loop space of $M$ whose critical points correspond to closed Reeb orbits. It is symmetric with respect to a suitable $L^2$-product. A periodic orbit $\gamma$ is non-degenerate if and only if $0$ does not lie in the spectrum of $\mathbf{A}_\gamma$, and more generally, if $\gamma$ is Morse--Bott and lies in a Morse--Bott submanifold $N$, then $\dim \ker \mathbf{A}_\gamma=\dim N - 1$. Under a choice of unitary trivialization $\tau$ of $\gamma^*\xi$, this operator looks like 
$$
\mathbf{A}_\gamma:W^{1,2}(S^1,\mathbb{R}^{2n})\rightarrow L^2(S^1,\mathbb{R}^{2n})
$$
$$
\mathbf{A}_\gamma =-i\partial_t-S(t),
$$
where $S$ is a smooth loop of symmetric matrices (the coordinate representation of $-TJ\nabla R$), which comes associated to a trivialization of $\mathbf{A}_\gamma$. When $\gamma$ is non-degenerate, its Conley--Zehnder index with respect to $\tau$ is defined to be the Conley--Zehnder index of the path of symplectic  matrices $\Psi(t)$ satisfying $\dot{\Psi}(t)=iS\Psi(t)$, $\Psi(0)=\mathds{1}$. We denote this by $\mu_{CZ}^\tau(\gamma)=\mu_{CZ}(\Psi)$.

We will consider, for cylindrical $J$, punctured $J$-holomorphic curves $u:(\dot{\Sigma},j)\rightarrow (\mathbb{R}\times M,J)$ in the symplectization of a stable Hamiltonian manifold $M$, where $\dot{\Sigma}=\Sigma \backslash \Gamma$, $(\Sigma,j)$ is a compact connected Riemann surface, and $u$ satisfies the nonlinear Cauchy--Riemann equation $du \circ j = J \circ u$. We will also assume that $u$ is asymptotically cylindrical, which means the following. Partition the punctures into \emph{positive} and \emph{negative} subsets $\Gamma = \Gamma^+\cup \Gamma^-$, and at each $z \in \Gamma^\pm$, choose a biholomorphic identification of a punctured neighborhood of $z$ with the half-cylinder $Z_\pm$, where $Z_+ = [0, \infty) \times S^1$ and $Z_- = (-\infty, 0] \times S^1$. Then writing $u$ near the puncture in cylindrical coordinates $(s, t)$, for $|s|$
sufficiently large, it satisfies an asymptotic formula of the form 
$$u \circ\phi(s, t) = exp_{(Ts,\gamma(T t))} h(s, t)$$
Here $T > 0$ is a constant, $\gamma : \mathbb{R} \rightarrow M$ is a $T$-periodic Reeb orbit, the exponential map is defined with respect to any $\mathbb{R}$-invariant metric on $\mathbb{R} \times M$, $h(s, t) \in \xi_{\gamma(T t)}$ goes to $0$ uniformly in $t$ as $s \rightarrow \pm \infty$ and $\phi : Z_\pm \rightarrow Z_\pm$ is a smooth embedding such that $\phi(s, t) - (s + s_0, t + t_0) \rightarrow 0$ as $s \rightarrow \pm \infty$ for some constants $s_0 \in \mathbb{R}$, $t_0 \in S^1$. We will refer to punctured asymptotically cylindrical $J$-holomorphic curves simply as $J$-holomorphic curves.

Observe that, for any closed Reeb orbit $\gamma$ and cylindrical $J$, the trivial cylinder over $\gamma$, defined as $\mathbb{R}\times \gamma$, is $J$-holomorphic. 

The Fredholm index of a punctured holomorphic curve $u$ which is asymptotic to non-degenerate Reeb orbits in a $(2n+2)$-dimensional symplectization $W^{2n+2}=\mathbb{R} \times M$ is given by the formula 

\begin{equation} \label{index} 
ind(u)=(n-2)\chi(\dot{\Sigma})+2c_1^\tau(u^*TW)+\mu^\tau_{CZ}(u)
\end{equation} 

Here, $\dot{\Sigma}$ is the domain of $u$, $\tau$ denotes a choice of trivializations for each of the bundles $\gamma_z^*\xi$, where $z \in \Gamma$, at which $u$ approximates the Reeb orbit $\gamma_z$. The term $c_1^{\tau}(u^*TW)$ is the relative first Chern number of the bundle $u^*TW$. In the case $W$ is $2$-dimensional, this is defined as the algebraic count of zeroes of a generic section of $u^*TW$ which is asymptotically constant with respect to $\tau$. For higher-rank bundles, one determines $c_1^\tau$ by imposing that $c_1^\tau$ is invariant under bundle isomorphisms, and satisfies the Whitney sum formula (see e.g.\ \cite{Wen8}). The term $\mu^\tau_{CZ}(u)$ is the total Conley--Zehnder index of $u$, given by 

$$
\mu^\tau_{CZ}(u)=\sum_{z \in \Gamma^+}\mu^\tau_{CZ}(\gamma_z) - \sum_{z \in \Gamma^-}\mu^\tau_{CZ}(\gamma_z)
$$ 

Given a $\mathcal{H}$-compatible $J$, and a $J$-holomorphic curve $u$ in $\mathbb{R}\times M$, the expression $u^*\Omega$ is a non-negative integrand, and one can define its $\Omega$-energy $$\mathbf{E}(u)=\int u^*\Omega$$ It is non-negative, and vanishes if and only if $u$ is a (multiple cover of) a trivial cylinder. 

\section{Algebraic torsion computations}\label{modelAsec}

\subsection{Construction of the model contact manifolds}
\label{model}

In this section, we construct the contact manifolds $M$ of Theorem \ref{thm5d}, making use of a cylindrical semi-filling $(Y\times I,d\alpha)$. We will use the ``double completion'' construction, originally appearing in \cite{LVHMW}. While very geometrically flavoured, this construction has the effect of endowing $M$ with a contact structure and an explicit deformation to a SHS, by viewing it as a contact-type hypersurface in a non-compact Liouville manifold. The contact form thus obtained will be degenerate, and a standard Morse function technique as in \cite{Bo} will be necessary.

\vspace{0.5cm}

\begin{figure}[t]\centering
\includegraphics[width=0.60\linewidth]{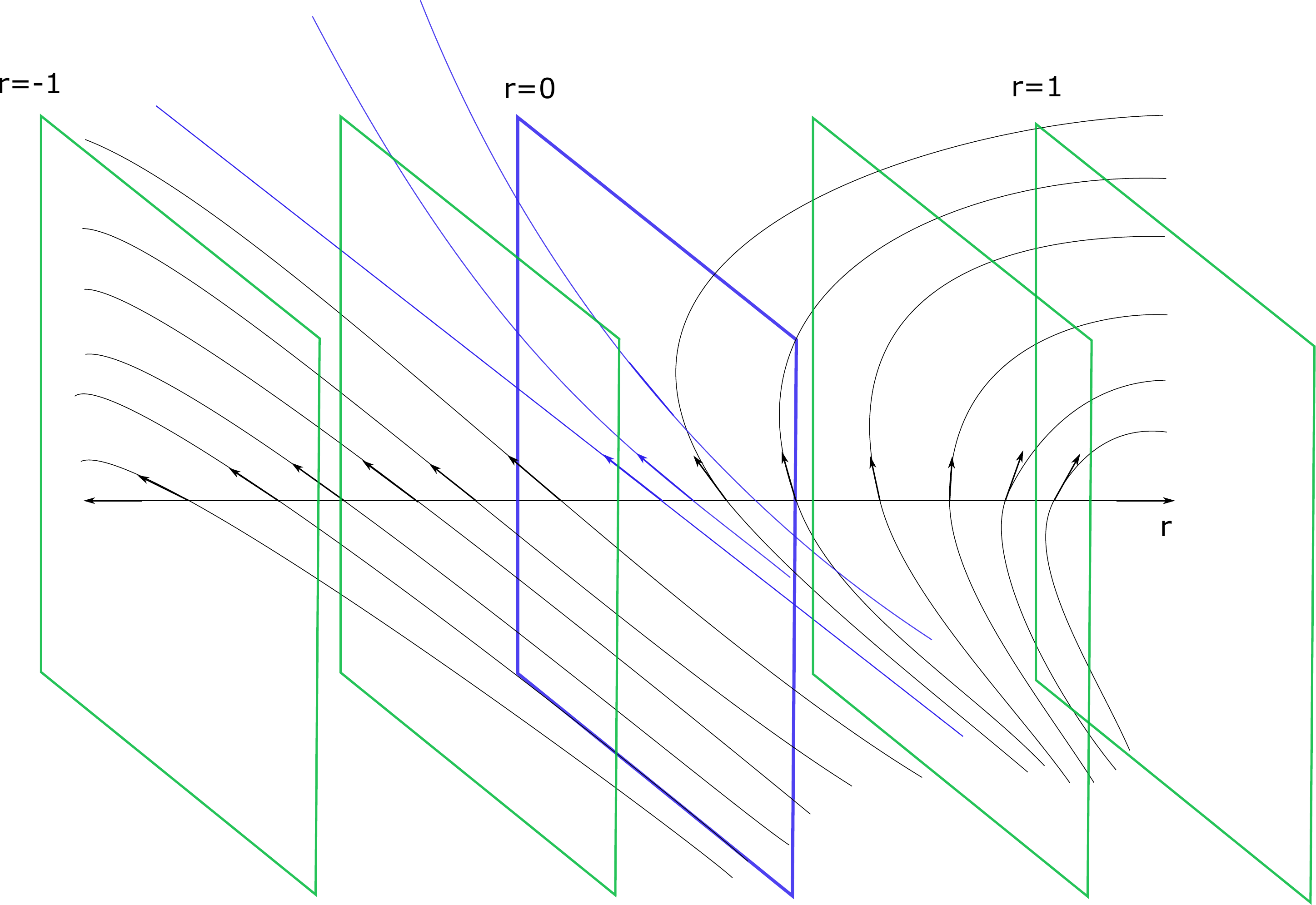}
\caption{\label{LouvilleDomain} The qualitative behaviour of the flow of the Liouville vector field $V$ on any cylindrical Liouville semi-filling, for which the central slice ($r=0$) is invariant. One may informally think of such a Liouville domain as being obtained by gluing two negative symplectizations along a ``non-contact hypersurface''.}
\end{figure}

Let $Y$ be a closed $(2n-1)$-manifold such that $(Y \times I,d\alpha)$ is a Liouville domain, for some exact symplectic form $d\alpha \in \Omega^2(Y \times I)$ (recall that throughout this paper, $I$ will denote the interval $[-1,1]$). See Figure \ref{LouvilleDomain} for a qualitative description.
We will assume that the Liouville form $\alpha=\{\alpha_r\}_{r \in I}$ is given by a 1-parameter family of $1$-forms in $Y$, which is the case for all known examples of cylindrical semi-fillings. In particular, we get that $\alpha(\partial_r)=0$. We can write the symplectic form as $$d\alpha=d\alpha_r + dr\wedge \frac{\partial \alpha_r}{\partial r}$$ The Liouville vector field $V$, defined to be $d\alpha$-dual to $\alpha$, points outwards at each boundary component, and hence, using its flow, we can choose our coordinate $r \in I$ so that $V$ agrees with $\pm \partial_r$ near the boundary $\partial(Y \times I)=Y \times \{\pm 1\}=:Y_\pm$. Therefore, we can assume that $\alpha= e^{\pm r-1}\alpha_\pm$ on $Y \times [-1,-1+\delta)$ and $Y \times (1-\delta,1]$, respectively, for some small $\delta>0$. Then $Y_\pm$ carries a contact structure $\xi_\pm=\ker\alpha_\pm$, where $\alpha_\pm = i_Vd\alpha|_{T(Y_\pm)}=\alpha|_{T(Y_\pm)}$. The behaviour of $V$ near the ends necessarily implies that there are values $r \in I$ such that $V|_{Y\times \{r\}}$ lies in $TY$, and hence $Y_r:=Y \times \{r\}$ is not a contact type hypersurface. The slices $Y_r$ which are of contact type inherit a contact structure $\xi_r=\ker \alpha_r$ and the resulting Reeb vector field $R_r$ satisfies $R_r=e^{1\mp r}R_\pm$ in the respective components of $\{|r| >1-\delta\}$, where $R_\pm$ is the Reeb vector field of $\alpha_\pm=\alpha_{\pm 1}$. We shall assume throughout that the only non-contact type slice is $Y_0$, so that $\alpha_r$ is a contact form for every $r \neq 0$. Also, we shall make the convention that whenever we deal with equations involving $\pm$'s and $\mp$'s, one has to interpret them as to having a different sign according to the region (the ``upper'' sign denotes the ``plus'' region, and the ``lower'', the ``minus'' region).

\vspace{0.5cm}

Let now $M= Y \times \Sigma$ be a product $(2n+1)$-manifold, where $\Sigma$ is the orientable genus $g$ surface obtained by gluing a connected genus $0$ surface with $k$ boundary components $\Sigma_-$, to a connected genus $g-k+1>0$ surface with $k$ boundary components $\Sigma_+$ along the boundary, by an orientation preserving map. The surface $\Sigma$ then inherits the orientation of $\Sigma_-$, which is opposite to the one in $\Sigma_+$. On each boundary component of $\partial \Sigma_\pm$, choose collar neighbourhoods $\mathcal{N}(\partial \Sigma_\pm)=(-\delta,0]\times S^1$ (for the same $\delta$ as before), and coordinates $(t_\pm,\theta_\pm) \in \mathcal{N}(\partial \Sigma_\pm)$, so that $\partial \Sigma_\pm=\{t_\pm = 0\}$. 

We will consider $\Sigma_-$ and $\Sigma_+$ to be attached at each of the $k$ boundary components by a cylinder $I\times S^1$, so that $M$ at this region is the disjoint union of $k$ copies of $Y \times I\times S^1$, with the $Y \times I$ identified with the Liouville domain above. We write the points of $M$ here as $(y,r,\theta)$, where the $\theta \in S^1$ coordinate can be chosen to coincide with $\theta_\pm$ where the gluing takes place. We shall therefore drop the subscript $\pm$ when talking about the $\theta$ coordinate. Denote also $$\mathcal{N}(-Y):=Y \times [-1,-1+\delta)\times S^1,\; \mathcal{N}(Y):=Y \times (1-\delta,1]\times S^1,$$ in the above identification. 

We have $$M=M_Y \cup M_P^{\pm},$$ where $M_Y=\bigsqcup^k Y \times I\times S^1$ is a region gluing $M_P^{\pm}=Y_\pm \times \Sigma_{\pm}$ together (recall Figure \ref{themanifold}). We shall refer to them as the \emph{spine} or \emph{cylindrical region}, and the \emph{positive/negative paper}, respectively. We have fibrations $$\pi_Y: M_Y \rightarrow Y \times I$$$$ \pi_P^\pm: M_P^{\pm} \rightarrow Y_\pm,$$ with fibers $S^1$ and $\Sigma_\pm$, respectively, and hence can be given the structure of a SOBD (see \cite{Mo2} for a definition). 

\vspace{0.5cm}

\begin{figure}[t]\centering \includegraphics[width=0.68\linewidth]{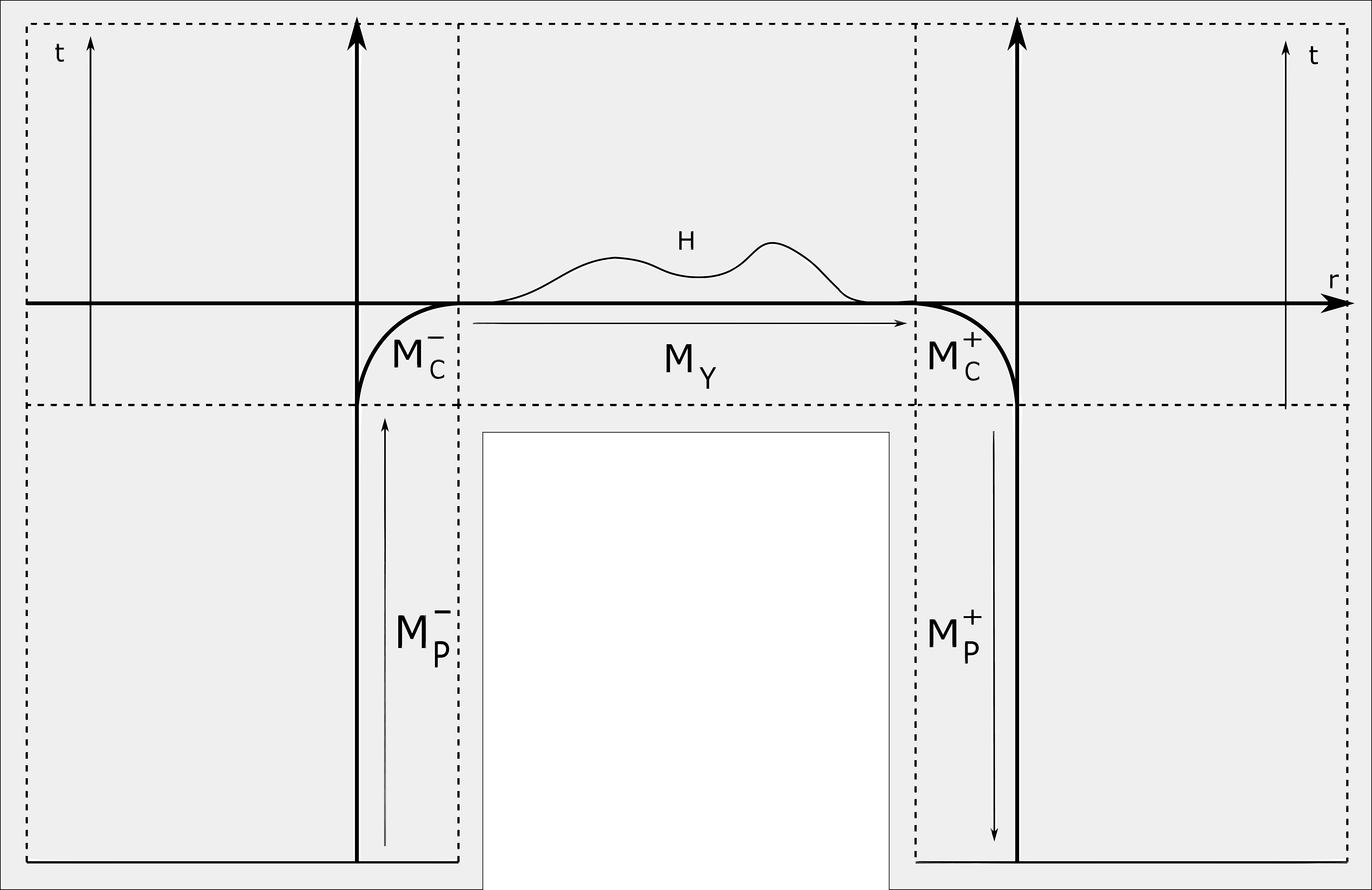}
\caption{\label{doublecompletion2} The double completion, and a Morse function $H$ along the spine.}
\end{figure}

We now construct an open manifold containing $M$ as a contact-type hypersurface. Denote by $\Sigma_\pm^\infty$ the open manifolds obtained from $\Sigma_\pm$ by attaching cylindrical ends of the form $(-\delta,+\infty) \times S^1$ at each boundary component, where the subset $(-\delta,0] \times S^1$ coincides with the collar neighbourhoods chosen above. The coordinates $t_\pm$ and $\theta$ extend to these ends in the obvious way, and we shall refer to the cylindrical ends as $\mathcal{N}(\partial \Sigma_\pm^\infty)$. We also consider the cylinder $\mathbb{R} \times S^1$ obtained by enlarging the cylindrical region $I\times S^1$ we had above. Denote then $$M_P^{\pm,\infty}=Y \times \Sigma_\pm^\infty$$$$M_Y^\infty= Y \times \mathbb{R} \times S^1$$$$ \mathcal{N}^\infty(-Y)=Y\times (-\infty, -1+\delta)\times S^1,\;\mathcal{N}^\infty(Y)=Y \times (1-\delta,+\infty) \times S^1$$ and define the \emph{double completion} of $E$ to be $$E^{\infty,\infty}=(-\infty,-1+\delta) \times M_P^{-,\infty} \;\bigsqcup\; (1-\delta,+\infty) \times M_
P^{+,\infty}\; \bigsqcup \; (-\delta,+\infty)\times M_Y^\infty /\sim,$$ where we identify $(r,y,t_-,\theta)\in (-\infty,-1+\delta)\times Y \times \mathcal{N}(\partial\Sigma_-^\infty)$ with $(t,y,r,\theta) \in (-\delta,+\infty) \times \mathcal{N}^\infty(-Y)$ if and only if $t=t_-$, and $(r,y,t_+,\theta)\in (1-\delta,+\infty)\times Y \times \mathcal{N}(\partial\Sigma_+^\infty)$ with $(t,y,r,\theta) \in (-\delta,+\infty) \times \mathcal{N}^\infty(Y)$ if and only if $t=t_+$  (see Figure \ref{doublecompletion2}). By definition, the $t$ coordinate coincides with the $t_\pm$ coordinates, where these are defined, so we shall again drop the $\pm$ subscripts from the variables $t_\pm$. Note also that the $r$ coordinate is globally defined, whereas $t$ is not. Denote then by $E^{\infty,\infty}(t)$ the region of $E^{\infty,\infty}$ where the coordinate $t$ is defined.

Choose now $\lambda_\pm$ to be Liouville forms on the Liouville domains $\Sigma^\infty_\pm$, such that $\lambda_\pm=e^td\theta$ on $\mathcal{N}(\partial \Sigma_\pm^\infty)$. This last expression makes sense in the region of $E^{\infty,\infty}$ where both $\theta$ and $t$ are defined, and where they are not, the form $\lambda_\pm$ makes sense. So this yields a globally defined 1-form $\lambda \in \Omega^1(E^{\infty,\infty})$ which coincides with $\lambda_\pm$ where these are defined. Also, the same argument works for $\alpha$, so that we get a global $\alpha \in \Omega^1(E^{\infty,\infty})$. 

For $K\gg 0$ a big constant, $\epsilon>0$ a small one, and $L\geq 1$, choose a smooth function $$\sigma=\sigma^{L}_{\epsilon,K}:\mathbb{R}\rightarrow \mathbb{R}^+$$ satisfying
\begin{itemize}
\item $\sigma \equiv K$ on $\mathbb{R}\backslash[-L,L]$.
\item $\sigma \equiv \epsilon$ on $[-L+\delta,L-\delta]$.
\item $\sigma^\prime(r)<0$, for $r\in (-L,-L+\delta)$, $\sigma^\prime(r) > 0$ for $r\in(L-\delta,L)$. 
\end{itemize}

We have that the 1-form $\sigma \alpha$ is Liouville on $Y \times \mathbb{R}$. Indeed, if $dvol$ is a positive volume form in $Y$ with respect to the $\alpha_-$-orientation, we may write 
$$
d\alpha^n=dvol\wedge dr,\; \alpha_r\wedge d\alpha_r^{n-1}=link(\alpha_r)dvol,
$$
where the last equation defines a \emph{self-linking} function $r \mapsto link(\alpha_r)$, whose sign is opposite to that of $r \in \mathbb{R}$. Then
$$
d(\sigma\alpha)^n=\sigma^{n-1}(\sigma - n \sigma^\prime link(\alpha_r))dvol\wedge dr
$$
Tracking the signs, one checks that the above expression is positive.

The associated Liouville vector field is  
\begin{equation} \label{Vsig}
V_\sigma:=\frac{\sigma}{\sigma + \sigma^\prime dr(V)}V=\left\{\begin{array}{ll} V,& \mbox{ on } (\mathbb{R}\backslash (-L,L))\cup [-L+\delta,L-\delta]\\
\frac{\sigma}{\sigma + \sigma^\prime}\partial_r,& \mbox{ on } (L-\delta,L)\\-\frac{\sigma}{\sigma - \sigma^\prime}\partial_r,& \mbox{ on } (-L,-L+\delta)
\end{array}\right.
\end{equation} 

Observe that $V_\sigma$ is everywhere \emph{positively} colinear with $V$. 

After extending the form $\sigma \alpha$ to $E^{\infty,\infty}$ in the natural way, one checks that  $$\lambda_\sigma:=\lambda_\sigma^{L}:=\sigma\alpha + \lambda$$ is a Liouville form on $E^{\infty,\infty}$. Denote $$\omega_\sigma:=\omega_\sigma^{L}:=d\lambda_\sigma=\sigma^\prime dr \wedge \alpha + \sigma d\alpha + d\lambda,$$ which is symplectic. Denote by $X_\sigma$ the associated Liouville vector field.

 If $X_\pm$ denotes the Liouville vector field on $\Sigma_\pm^\infty$ which is $d\lambda_\pm$-dual to $\lambda_\pm$, coinciding with $\partial_t$ in $\mathcal{N}(\partial \Sigma_\pm^\infty)$, we can define a smooth vector field on $E^{\infty,\infty}$ by 
\begin{equation}
\label{X}
X=\left\{\begin{array}{ll}
  X_+,& \mbox{ on } \{r>1-\delta\} \\
  X_-, &\mbox{ on } \{r<-1+\delta\} \\
  \partial_t,& \mbox{ on } \{|r|<1\}\\
  \end{array}\right.
\end{equation}
Then $$X_\sigma= X+V_\sigma$$   

Denote $$E^{L,Q}=E^{\infty,\infty}/(\{|r|>L\} \cup \{t > Q\}),$$ for $Q \geq 0$, and $L \geq 1$. We have its ``horizontal'' and ``vertical'' boundaries $$\widetilde{M}_Y^{L,Q}:=\partial_hE^{L,Q} := \{t=Q\} \cap \{r \in [-L,L]\}$$$$\widetilde{M}_P^{L,Q}:=\partial_vE^{L,Q} := \widetilde{M}_P^{-,L,Q} \bigsqcup \widetilde{M}_P^{+,L,Q},$$ where $$\widetilde{M}_P^{-,L,Q} :=\{r=-L\}\cap \{t \leq Q\}$$$$\widetilde{M}_P^{+,L,Q}:= \{r=L\}\cap \{t \leq Q\}$$
The manifold $$\widetilde{M}^{L,Q}:=\partial E^{L,Q}:= \partial_h E^{L,Q} \cup \partial_v E^{L,Q},$$ is then a manifold with corners $$\partial_h E^{L,Q} \cap \partial_v E^{L,Q}=\{|r|=L\}\cap \{t=Q\}$$  
  
One has $$X_\sigma=\pm\frac{\sigma}{\sigma \pm \sigma^\prime}\partial_r + \partial_t$$ in the corresponding components of the region $\{|r|>L-\delta\}\cap \{t\geq -\delta\}$. This means that $X_\sigma$ will be transverse to the smoothening of $\partial E^{L,Q}$ that we shall now construct.

\vspace{0.5cm}

\begin{figure}[t] \centering \includegraphics[width=0.60\linewidth]{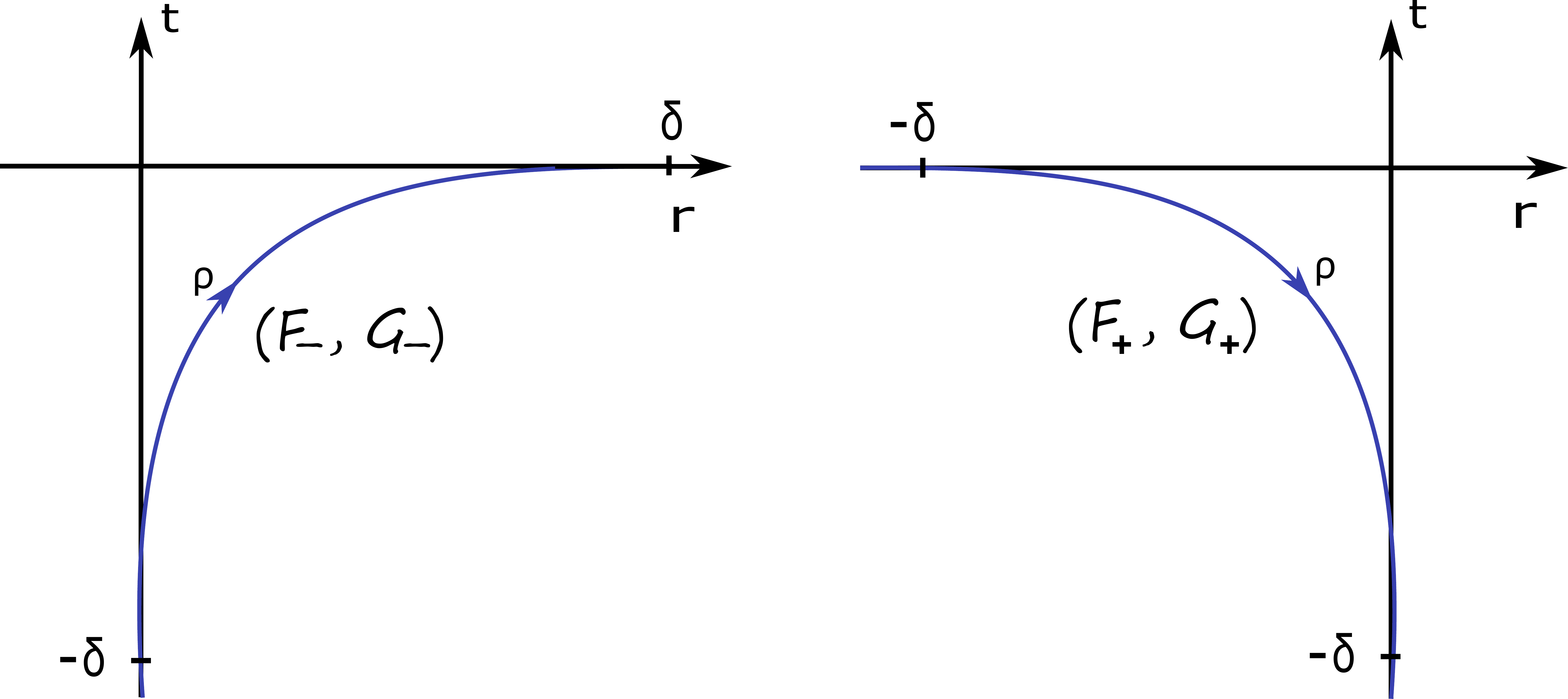}
\caption{\label{SmoothendCorner} The paths $\rho \mapsto (F_\pm(\rho),G_\pm(\rho))$.}
\end{figure}

Choose smooth functions $F_\pm,G_\pm:(-\delta,\delta)\rightarrow (-\delta,0]$ such that $$\left\{\begin{array}{ll}
 (F_+(\rho),G_+(\rho))=(\rho,0),\;(F_-(\rho),G_-(\rho))=(0,\rho),& \mbox{ for } \rho \leq -\delta/3\\
 G_+^\prime(\rho)<0, F_-^\prime(\rho)>0,& \mbox{ for } \rho > -\delta/3\\
 G_-^\prime(\rho)>0, F_+^\prime(\rho)>0,& \mbox{ for } \rho < \delta/3\\
 (F_+(\rho),G_+(\rho))=(0,-\rho),\;(F_-(\rho),G_-(\rho))=(\rho,0),& \mbox{ for } \rho \geq \delta/3\\
\end{array}\right.$$ See Figure \ref{SmoothendCorner}.

We now smoothen the corner $\partial_h E^{L,Q} \cap \partial_v E^{L,Q}$ by substituting the region $$\partial E^{L,Q} \cap \left(\{t \in (Q-\delta,Q]\}\cup \{|r|\in (L-\delta,L]\}\right),$$ which contains the corners, with the smooth manifold $$M_{C}^{\pm,L,Q}:=\{(r=\pm L+F_\pm(\rho), t = Q + G_\pm(\rho),y,\theta): (\rho,y,\theta) \in (-\delta,\delta) \times Y \times S^1\}$$

The smoothened boundary can then be written as $$M^{L,Q}:= M_P^{\pm,L,Q} \bigcup M_C^{\pm,L,Q} \bigcup M_Y^{L,Q},$$ where $$M_P^{\pm,L,Q} = \widetilde{M}_P^{\pm,L,Q} \cap \{t<Q-\delta/3\}$$$$M_Y^{L,Q}= \widetilde{M}_Y^{L,Q} \cap \{|r|<L-\delta/3\}$$ The Liouville vector field $X_\sigma$ is transverse to this manifold, so that we get a contact structure on $M^{L,Q}$ given by $$\xi^{L,Q}= \ker(\lambda^{L}_\sigma|_{TM^{L,Q}})$$ Observe that $M^{L,Q}$ is canonically diffeomorphic to $M$. So, this actually yields a contact structure on $M$. By construction, we have non-empty intersections $$M_P^{\pm,L,Q}\cap M_C^{\pm,L,Q} = \{r=\pm L\}\cap \{t \in (Q-\delta,Q-\delta/3)\}$$$$ M_Y^{L,Q}\cap M_C^{\pm,L,Q}=\{t=Q\}\cap \{|r| \in (L-\delta,L-\delta/3)\}$$

We shall construct a stable Hamiltonian structure on $M^{L,Q}$ which arises as a deformation of the above contact structure, such that both coincide on $M_Y^{L,Q}$, as follows.

Choose a smooth function $\beta: \mathbb{R}\rightarrow [0,1]$ such that $\beta(t)=0$ for $t \leq -\delta+\delta/9$, $\beta(t)=1$ for $t\geq -2\delta/3-\delta/9$, and $\beta^\prime \geq 0$. Set 
\begin{equation}\label{Z}
Z=\left\{\begin{array}{ll} V_\sigma+\beta(t)X,& \mbox{ in the region } E^{\infty,\infty}(t) \mbox{ where } t \mbox{ is defined}\\ 
\frac{\sigma}{\sigma+\sigma^\prime}\partial_r,& \mbox{ in } E^{\infty,\infty}(t)^c \cap \{r > 1-\delta\}\\
-\frac{\sigma}{\sigma-\sigma^\prime}\partial_r,& \mbox{ in } E^{\infty,\infty}(t)^c \cap \{r < -1+\delta\},
\end{array}\right. 
\end{equation}
which yields a smooth vector field on $E^{\infty,\infty}$, a deformation of $X_\sigma$. Then $Z$ is still transverse to $M^{L,Q}$ and is stabilizing, so that the pair $$\mathcal{H}:=\mathcal{H}^{L,Q}:=(\Lambda^{L,Q}=i_Z\omega_\sigma|_{TM^{L,Q}},\Omega^{L,Q}=\omega_\sigma|_{TM^{L,Q}})$$ yields a stable Hamiltonian structure on $M^{L,Q}$. For $Q=0$, $(M^{L,0}_Y, \Lambda^{L,0})$ can be seen as the \emph{contactization} of the Liouville domain $(Y\times I,\epsilon d\alpha)$.  

\vspace{0.5cm}

Along $M_Y^{L,Q}$ the Reeb vector field is given by $R^{L,Q}= \frac{\partial_\theta}{e^Q},$ which is degenerate, and the space of Reeb orbits is identified with $Y \times [-L,L]$. We consider two perturbation approaches: Morse, and Morse-Bott. In the first approach we choose $H_L: Y \times [-L,L] \rightarrow \mathbb{R}^{\geq 0}$ to be Morse, depending only in $r$ near $r=\pm L$, satisfying $\partial_rH_L\leq 0$ near $r=L$, $\partial_rH_L\geq 0$ near $r=-L$, and vanishing as one approaches $r=\pm L$. In the second approach, we choose $H_L$ to depend only on $r$ \emph{globally}, with respect to which it is a Morse function. 

\vspace{0.5cm}

If $t\mapsto \Phi^t_Z$ denotes the flow of $Z$, choose $\epsilon > 0$ sufficiently small so that the manifold $$M^{\epsilon; L,Q}:=\{\Phi_Z^{\epsilon H_L(x)}(x) \in E^{\infty,\infty}: x \in M^{L,Q}\}$$ is still transverse to $Z$. We have a stable Hamiltonian structure $$\mathcal{H}_\epsilon:=\mathcal{H}_\epsilon^{L,Q}:=(\Lambda_\epsilon^{L,Q},\Omega_\epsilon^{L,Q}):=(i_Z\omega_\sigma|_{TM^{\epsilon;L,Q}},\omega_\sigma|_{TM^{\epsilon;L,Q}}),$$ and a decomposition $$M^{\epsilon;L,Q}=M_Y^{\epsilon;L,Q} \cup M_C^{\epsilon;\pm,L,Q} \cup M_P^{\epsilon;\pm,L,Q},$$ where each component is the perturbation of the corresponding component of $M^{L,Q}$. 

Along the region $M_C^{\epsilon;\pm,L,Q}$ the new coordinates are  $$r=F_\pm^{\epsilon;L}(\rho)=\Phi_{V_\sigma}(\epsilon H_L(\cdot,\pm L+F_\pm(\rho)),\pm L+F_\pm(\rho))$$$$ t=G_\pm^{\epsilon;L,Q}(\rho)=Q+G_\pm(\rho)+ \epsilon H_L(\cdot,\pm L+F_\pm(\rho)),$$ where $\Phi_{V_\sigma}(s,\cdot)$ is the time $s$ flow of $V_\sigma$.

We then have
$$
\Lambda_\epsilon^{L,Q}=\left\{
\begin{array}{ll}\Lambda^{L,Q}=K \alpha_\pm, & \mbox{ in } M_P^{\epsilon;\pm,L,Q}=M_P^{\pm,L,Q}\\
e^{\epsilon H_L}(\epsilon\alpha + e^Qd\theta), & \mbox{ in } M_Y^{\epsilon;,L,Q}\\
\sigma(F_\pm^{\epsilon;L}(\rho)) e^{\pm F_\pm^{\epsilon;L}-L} \alpha_\pm + \beta(G_\pm^{\epsilon;L,Q}(\rho))e^{G_\pm^{\epsilon;L,Q}(\rho)}d\theta,& \mbox{ in } M_C^{\epsilon;\pm,L,Q}
\end{array}\right.
$$
One can similarly write down $\Omega_\epsilon^{L,Q}$ explicitly. 

The Reeb vector field $R^{L,Q}_\epsilon$ associated to this stable Hamiltonian structure is 

$$R_\epsilon^{L,Q}=\left\{\begin{array}{ll} R^{L,Q}=\frac{R_\pm}{K}, & \mbox{ in } M_P^{\epsilon;\pm,L,Q}\\
e^{-\epsilon H_L-Q}\left((1+\epsilon \alpha(X_{H_L}))\partial_\theta-e^QX_{H_L}\right),& \mbox{ in } M_Y^{\epsilon;L,Q}\\
\frac{1}{\Phi^{\epsilon; L,Q}_\pm}\left(e^{\mp F_\pm^{\epsilon;L}+L}(G_\pm^{\epsilon;L,Q})^\prime R_\pm- e ^{-G^{\epsilon;L,Q}_\pm}(\sigma^\prime\pm \sigma)(F_\pm^{\epsilon;L})(F_\pm^{\epsilon;L})^\prime\partial_\theta\right), & \mbox{ in } M_C^{\epsilon;\pm,L,Q}   
\end{array}\right.
$$ 
where $X_{H_L}$ is the Hamiltonian vector field on $Y \times I$ associated to $H_L$, defined by $i_{X_{H_L}}d\alpha = -dH_L$, and \begin{equation} \label{Phieps} 
\Phi^{\epsilon; L,Q}_\pm(\rho)= \sigma(F_\pm^{\epsilon;L}(\rho)) (G_\pm^{\epsilon;L,Q})^\prime(\rho) - \beta(G_\pm^{\epsilon;L,Q}(\rho))(\sigma^\prime \pm \sigma)(F_\pm^{\epsilon;L}(\rho))(F_\pm^{\epsilon;L})^\prime(\rho)
\end{equation}

One can check that $\Phi^{\epsilon; L,Q}_\pm$ has sign which is \emph{opposite} to its subscript. Observe that critical points $(y,r)$ of $H$ give rise to closed Reeb orbits of the form $\gamma_p:=\{p\} \times S^1 \subseteq \mbox{crit}(H_L)\times S^1$. If we are taking the Morse approach, we have only a finite number of such orbits, and they are non-degenerate. Choosing $H_L$ to be $C^1$-small has the effect of making the vector field $X_{H_L}$ also small, so that the closed orbits which do not arise from critical points of $H_L$ have large period, including the ones not contained in $M_Y^{L,Q}$. So, taking any large (but fixed) $T\gg 0$, we can choose $H_L$ small enough so that all the periodic Reeb orbits up to period $T$ are of the form $\gamma_p^l$, for $p \in \mbox{crit}(H_L)$, and $l\leq N$, for some covering threshold $N$ depending on $T$. For the Morse--Bott case, we obtain $Y$-families of Morse-Bott orbits for each critical point of $H_L$. 

\vspace{0.5cm}

\begin{remark} \label{remark11}
\begin{enumerate}[wide, labelwidth=!, labelindent=0pt]
$\;$
\item One can check that $\lambda_\sigma(R_\epsilon)=1$ (recall that $\lambda_\sigma\vert_{M^{\epsilon;L,Q}}$ is the primitive of $\Omega_\epsilon=\omega_\sigma\vert_{M^{\epsilon;L,Q}}$). Therefore, for compatible almost complex structure $J$ and an asymptotically cylindrical $J$-holomorphic curve $u$ with positive/negative punctures $\Gamma^\pm$, the $\Omega_\epsilon$-energy of $u$ is
\begin{equation}
\label{action}
\int u^*\Omega_\epsilon=\int u^*d\lambda_\sigma=\sum_{z \in \Gamma^+}T_z-\sum_{z \in \Gamma^-}T_z,
\end{equation}
where $T_z$ is the action of the Reeb orbit corresponding to the puncture $z$. In particular, if the positive punctures correspond to critical points of $H_L$, then so will the negative ones.

\item By inspecting the expressions of the Reeb vector field we see that there are no contractible closed Reeb orbits for the SHS, if we assume this same condition for $R_\pm$. Moreover, the direction of the Reeb vector field does not change after perturbing back to sufficiently close contact data (cf.\ Section \ref{SHStoND} below), and so this also holds for the latter data. It follows that the isotopy class defined by the resulting contact structure is hypertight, and this shows the hypertightness condition of Theorem \ref{thm5d}.

\end{enumerate}
\end{remark}

\subsection{Compatible almost complex structure}
\label{ACS}

\paragraph*{Construction}

We set $L=1$, and $Q=0$, and drop the superscripts $L$ and $Q$ from all of the notation. We now define a suitable, though non-generic, almost complex structure $J=J_\epsilon$ on the symplectization $W^\epsilon=\mathbb{R} \times M^\epsilon$, where $$M^\epsilon=M^{\epsilon;1,0}=M^\epsilon_Y \bigcup M_P^{\epsilon;\pm} \bigcup M^{\epsilon;\pm}_C$$ It will be compatible with the stable Hamiltonian structure $\mathcal{H}_\epsilon$, and the fibers $\Sigma_\pm$ of our fibration $\pi_P^\pm$, the ``pages'', will lift as holomorphic curves. We will blur the distinction between $M=M^{0;1,0}$ and its diffeomorphic perturbed copy $M^\epsilon$ (as well as for $W$ and $W^\epsilon$), so that we are actually working on a fixed $M$ with a SHS which depends on $\epsilon$. 

\vspace{0.5cm}

Denote by $\xi_\epsilon:= \ker \Lambda_\epsilon$. We will define $J$ on $\xi_\epsilon$, in an $\mathbb{R}$-invariant way, and then simply set $J(\partial_a):=R_\epsilon$.

Choose a $d\alpha$-compatible almost complex structure $J_0$ on $Y\times I$, which is cylindrical in the cylindrical ends of $Y\times I$, so that, along these, it coincides with a $d\alpha_\pm$-compatible almost complex structure $J^\pm$ on $\xi_\pm$, and maps the Liouville vector field $V$ to $R_\pm$. Observe that the vector field $\partial_\theta$ is transverse to $\xi_\epsilon$ along $M_Y$. Therefore, we may then define $J_Y=\pi_Y^* J_0$ on $\xi_\epsilon\vert_{M_Y}$.

\vspace{0.5cm}

 Along the regions $\Sigma_\pm/\mathcal{N}(\partial\Sigma_\pm) \times Y\subseteq M_P^{\pm}$, and $\{t \in (-\delta,-\delta/3)\}\times Y\subseteq M_P^{\pm}$, the restriction of the projection $\pi_P^\pm: M_P^\pm\rightarrow \Sigma_\pm$ induces an isomorphism $d\pi_P^\pm: \xi_\epsilon/\xi_\pm \stackrel{\simeq}{\longrightarrow} T\Sigma_\pm$. Choose $j_\pm$ to be a $d\lambda$-compatible almost complex structures on $\Sigma_\pm$, so that $j_\pm(\partial_t)= K \partial_\theta$ in $\mathcal{N}(\partial \Sigma_\pm)$. Define $$J=J^\pm \oplus (\pi_P^\pm)^*j_\pm$$ on $\xi_{\epsilon}\vert_{M_P^\pm}$.

\vspace{0.5cm}

In $M_C^{\pm}$, we have  
$$\xi_\epsilon=\xi_\pm \oplus \langle v_1,v_2^\pm\rangle,$$ where 

\begin{equation}\label{wpm}
v_1=\partial_\rho,\;v_2^\pm=a_\pm(\rho)R_\pm + b_\pm(\rho) \partial_\theta,
\end{equation}
Here, 
$$
a_\pm=-\frac{\beta(G_\pm^\epsilon)}{e^{\pm F^\epsilon_\pm-1}\Phi_\pm^\epsilon}
$$ 
$$
b_\pm=\frac{\sigma(F_\pm^\epsilon)}{ e^{G_\pm^\epsilon}\Phi_\pm^\epsilon}
$$ 
where $\Phi_\pm^\epsilon$ is defined in (\ref{Phieps}). In the overlaps $M_C^{\pm}\cap M_Y$, one computes that $J(v_1)= g_\pm^Y(\rho) v_2^\pm$
where $g_\pm^Y:=\pm \frac{e^{\mp F^\epsilon\pm(\rho)+1}}{a_\pm(\rho)}=\mp\frac{\Phi_\pm^\epsilon}{\beta(G_\pm^\epsilon)}$, which is always positive. Similarly, in $M_P^{\pm}\cap M_C^{\pm}$, we have $J(v_1)=e^{\mp \rho} v_2^\pm$. Since $g_\pm^Y$ and $e^{\mp\rho}$ are both positive, we can now take any smooth positive functions $h_\pm:(-\delta,\delta) \rightarrow \mathbb{R}^+$ which coincide with $e^{\mp \rho}$ near $\rho=\pm \delta$ and with $g_\pm^Y$ near $\rho=\mp\delta$. We glue the two definitions by setting 
\begin{equation}
\label{w}
J(v_1)=h_\pm(\rho)v_2^\pm:=w^\pm,
\end{equation}
and we make $J$ agree with $J^\pm$ on $\xi_\pm$.

This gives a well-defined cylindrical $J$ in $\mathbb{R}\times M$.

\paragraph*{Compatibility.} One can check that $J$ is $\mathcal{H}_\epsilon$-compatible by straightforward computations \cite{Mo2}.

\begin{remark}\label{dlambda}$\;$
We observe that over $\mathbb{R}\times M_P^\pm$, where $\Lambda_\epsilon=K\alpha_\pm$, we have $d\Lambda_\epsilon(v,Jv)\geq 0$ , with equality if and only if $v^{\xi_\pm}=0$, so that the projection to $TY$ of $v$ lies in the span of $R_\pm$.
\end{remark}

\subsection{Finite energy foliation}
\label{Curves}

We will now consider the symplectization of our stable Hamiltonian manifold $(M, \mathcal{H}_\epsilon)$, given by $(W=\mathbb{R}\times M,\omega^\varphi_\epsilon=d(\varphi(a)\Lambda_\epsilon)+\Omega_\epsilon)$, where $\varphi \in \mathcal{P}$. We will construct a finite-energy foliation of $W$ by $J$-holomorphic curves, consisting of three distinct types, which we describe in the next theorem. This is an adaptation of the construction in \cite{Wen5}.

\begin{thm}
There exists a finite energy foliation of the symplectization $(W,\omega_\epsilon^\varphi)$ by simple $J$-holomorphic curves of the following types:

\begin{itemize}
\item \textbf{trivial cylinders $C_p$}, corresponding to Reeb orbits of the form $\gamma_p=\{p\} \times S^1$ for $p \in$ crit$(H)$, and which may be parametrized by 

$$C_p: \mathbb{R}\times S^1 \rightarrow \mathbb{R} \times M_Y =\mathbb{R} \times (Y \times I) \times S^1 $$$$C_p(s,t)=(s,p,e^{-\epsilon H(p)}t)$$ 

\item \textbf{flow-line cylinders $u^a_\gamma$}, parametrized by $$u^a_\gamma: \mathbb{R}\times S^1 \rightarrow \mathbb{R} \times M_Y =\mathbb{R} \times (Y \times I) \times S^1$$
\begin{equation}\label{trivcyl}
u^a_\gamma(s,t)=(a(s),\gamma(s),\theta(s)+t)
\end{equation}
for a proper function $a: \mathbb{R}\rightarrow \mathbb{R}$, a function $\theta:\mathbb{R}\rightarrow S^1$, and a map $\gamma: \mathbb{R}\rightarrow Y\times I$ satisfying 

\begin{equation}
\label{floweqs}
\begin{array}{l} 
\dot{\gamma}= \nabla H(\gamma(s))\\
\dot{a}(s)= e^{\epsilon H(\gamma(s))},\; a(0)=a\\
\dot{\theta}(s)=-\epsilon\alpha(\nabla H (\gamma(s))),\; \theta(0)=0
\end{array}
\end{equation}

Here, the gradient is computed with respect to the metric $g_{d\alpha,J_0}:= d\alpha(\cdot, J_0\cdot)$.

They have for positive/negative asymptotics the Reeb orbits corresponding to $$p_\pm:=\lim_{s\rightarrow \pm \infty} \gamma(s) \in \mbox{crit}(H)$$ 

\item \textbf{positive/negative page-like holomorphic curves $u^\pm_{y,a}$}, which consist of a trivial lift at symplectization level $a$ of a page $P_y^\pm:= \left(\Sigma_\pm\backslash\mathcal{N}(\partial\Sigma_\pm)\right) \times \{y\}$, for $y \in Y$, glued to $k$ cylindrical ends which lift the smoothened corners and then enter the symplectization of $M_Y$, asymptotically becoming a flow-line cylinder. They have $k$ positive asymptotics at Reeb orbits of the form $\gamma_p$, exactly one for each component of $M_Y$. The positive curves have genus $g-k+1>0$ and $k$ punctures, whereas the negative curves have genus $0$, and also $k$ punctures.  
\end{itemize}
\end{thm}

\begin{figure}[t]
\centering
\includegraphics[width=0.7\linewidth]{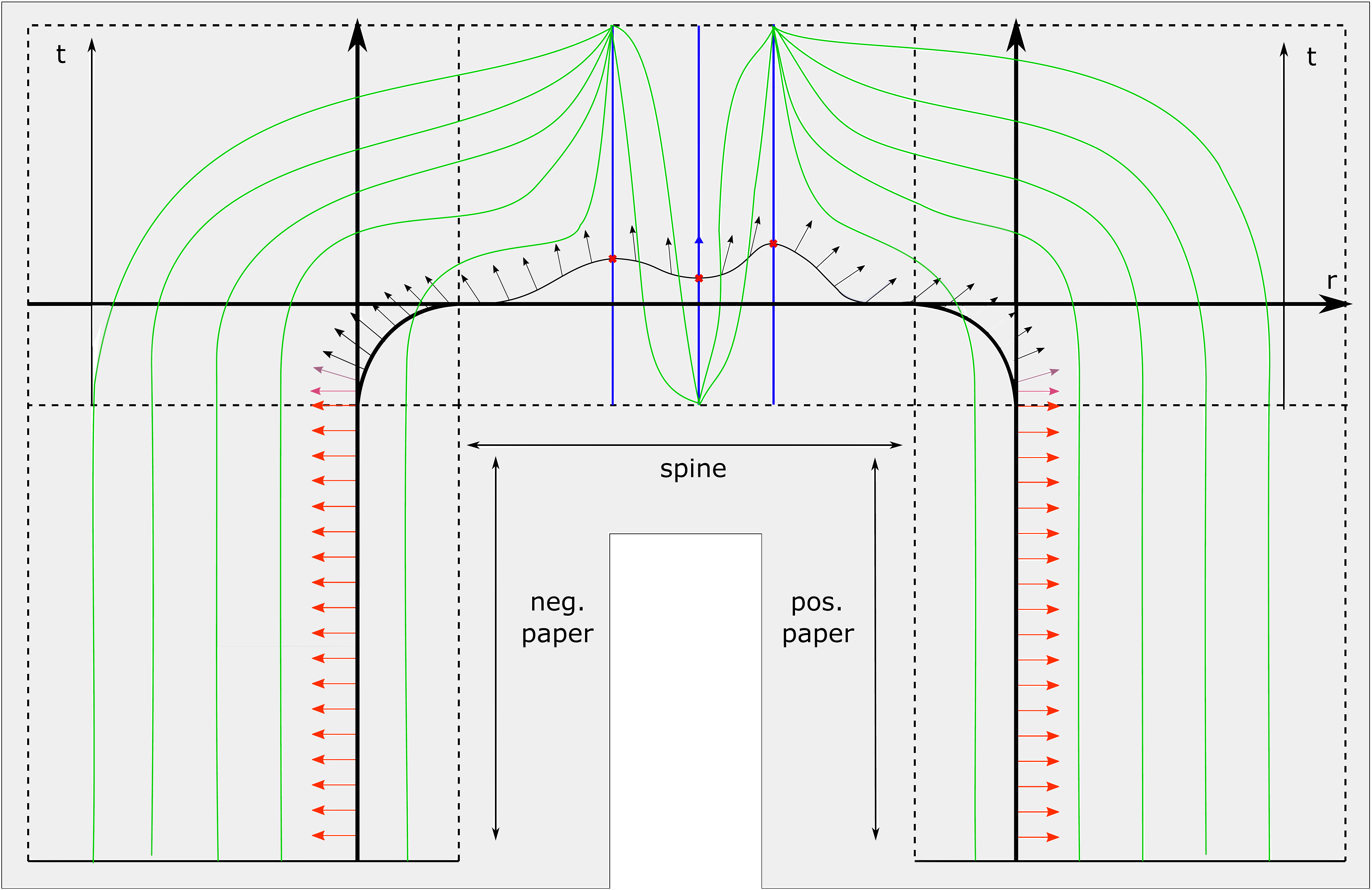}
\caption{\label{DoubleCompletion} The double completion $E^{\infty,\infty}$, containing $M$ and its perturbed version $M^\epsilon$ as contact type hypersurfaces. The foliation by holomorphic curves is shown in green (for the non-trivial curves) and blue (for the trivial cylinders).}
\end{figure}

\begin{remark}\label{gradrk} In the Morse-Bott case, one can show that $\alpha(\nabla H)=0$ for a suitable metric on $Y\times I$ (\cite{Mo2}, Remark 2.8), so that the function $\theta$ vanishes identically.
\end{remark}

Figure \ref{DoubleCompletion} summarizes the situation. We shall not distinguish the curves $u_{y,a}^{\pm}$ and $u^a_\gamma$ from the simple holomorphic curves that they parametrize, and we will drop the $a$ for the notation whenever we wish to refer to the equivalence class of the curves under $\mathbb{R}$-translation. A short computation shows that the flow-line cylinders are indeed holomorphic (see e.g.\ \cite{Mo2, Sie}). We now construct the page-like curves.

The pages $P^\pm_y:=(\Sigma_\pm\backslash\mathcal{N}(\partial \Sigma_\pm)) \times \{y\}$, $y \in Y$, clearly lift to a holomorphic foliation of the region $\mathbb{R} \times M_P^{\pm} \subseteq W$, which takes the form $\{\{a\} \times P^\pm_y :a \in \mathbb{R}, y \in Y\}$. We now glue cylindrical ends to these lifts.

\vspace{0.5cm}

We have that $Jv_1=h_\pm v_2^\pm$ and $R_\epsilon$ are both linear combinations of the vector fields $R_\pm$ and $\partial_\theta$ along $M_C^\pm$, with the coefficients only depending on $\rho$. Since these are not colinear, we have smooth functions $B,C:(-\delta,\delta)\rightarrow \mathbb{R}$ such that $$\partial_\theta = B R_\epsilon + CJv_1$$ One can in fact compute that the following expressions hold: 
\begin{equation}\label{BC}
\begin{array}{l}
B=e^{G_\pm^\epsilon}\beta(G_\pm^\epsilon)\\
C=\frac{e^{G_\pm^\epsilon}(G_\pm^\epsilon)^\prime}{h_\pm}
\end{array}
\end{equation}

We have that 
\begin{equation}\label{holcorner}
J\partial_\theta = -B\partial_a  -Cv_1=-B\partial_a-C\partial_\rho
\end{equation}
We conclude that 
$$\langle \partial_\theta , B\partial_a+C\partial_\rho\rangle = \langle \partial_\theta ,J\partial_\theta \rangle
$$

It follows that the distribution above has integral submanifolds which are unparametrized holomorphic curves. We can actually find holomorphic parametrizations given by 

$$w^\pm_{y,a}: (-\delta,\delta) \times S^1 \rightarrow  \mathbb{R} \times M_C^{\pm} =\mathbb{R}\times  Y \times (-\delta,\delta)\times S^1$$$$ (s,t)\mapsto (b(s),y,\rho(s),t),$$ 

for some fixed $y \in Y$, and functions $b,\rho:(-\delta,\delta)\rightarrow \mathbb{R}$ satisfying 

$$\dot{\rho}(s)=C(\rho(s)),\; \rho(-\delta)=-\delta$$
$$\dot{b}(s)=B(\rho(s)),\;b(-\delta)=a$$

The curve $w^\pm_{y,a}$ is indeed holomorphic. 

\vspace{0.5cm}

The curves $w^\pm_{y,a}$ glue with curves $u^a_\gamma$, which look like $u^a_\gamma(s,t)=(a(s),y(s),r(s),t)$ for some $y(s) \in Y$ such that $y(s)\equiv y$ near $r=\pm 1$ and $\lim_{s\rightarrow +\infty}y(s)=y^+$, and some $r:\mathbb{R}\rightarrow I$ with $\lim_{s\rightarrow +\infty}r(s)=r^+$, so that $(y^+,r^+)\in crit(H)$. We may then define a $J$-holomorphic curve $$u^\pm_{y,a}:=P^\pm_{y,a} \bigcup w^\pm_{y,a} \bigcup u^a_\gamma$$ which asympote $k$ Reeb orbits $\gamma_{y,a;i}^\pm=\{p_{y,a;i}^\pm\} \times S^1$, where $p_{y,a;i}^\pm \in crit(H)$, for $i=1,\dots,k$, and which have genus $g(u^-_{y,a})=0$ and $g(u^+_{y,a})=g-k+1$.  

\subsection{Index computations}
\label{indexcomp}

In this section, we compute the Fredholm index of the curves in the foliation.

\begin{thm} $\;$

\begin{enumerate}[wide, labelwidth=!, labelindent=0pt] 
\item After a sufficiently small Morse perturbation making Reeb orbits along $M_Y$ non-degenerate, we can find a natural trivialization $\tau$ of the contact structure along $\gamma_p$ (inducing a trivialization $\tau^l$ along all of its covers $\gamma_p^l$), and $N\in \mathbb{N}$, which depends on $H$ and grows as $H$ gets smaller, such that the Conley--Zehnder index of $\gamma_p^l$ is given by $$\mu^{\tau^l}_{CZ}(\gamma^l_p)=\mbox{ind}_p(H)-n,$$ for $l\leq N$.
  
\item In the Morse approach, the Fredholm indexes of the curves in our finite energy foliation are given by 
\begin{equation}
\begin{split}
\mbox{ind}(u_{y,a}^-)&=2n(1-k)+ \sum_{i=1}^k \mbox{ind}_{p^-_{y,a;i}}(H)\\ 
\mbox{ind}(u_{y,a}^+)&=2n(1-g-k)+\sum_{i=1}^k ind_{p_{y,a;i}^+}(H)\\
\mbox{ind}(u_\gamma^a)&=\mbox{ind}_{p_+}(H)-\mbox{ind}_{p_-}(H)
\end{split}
\end{equation}
\end{enumerate}
\end{thm} 

\begin{proof}See \cite{Mo2}. \end{proof}

\begin{remark} 
Since $ind_{p^+_{y,a;i}}(H)\leq 2n$ for every $i$, then $ind(u_{y,a}^+)\leq 2n(1-g)\leq 0$, since $g\geq 1$. This means these curves cannot possibly achieve transversality, and, after a perturbation making $J$ generic, they will disappear.
\end{remark}

\subsection{Fredholm regularity}
\label{Regularity}

In this section, we shall prove that the curves we have constructed are Fredholm regular.

In the Morse case, regularity of unbranched covers of flow-line cylinders can be reduced to the Morse--Smale condition for $H$. This fact is known to experts, and we shall omit the details (see \cite{Mo2}).     

For regularity of the other curves, we will assume the Morse--Bott situation, and prove regularity of the genus zero curves in the foliation. We will use the fact from \cite{Wen1} that Fredholm regularity is equivalent to the surjectivity of the normal component of the linearized Cauchy-Riemann operator, which is again a Fredholm operator. After some assumptions on our choice of coordinates and the Morse-Bott function $H$ (which we can always assume hold), and a suitable choice of normal bundle, we will explicitly write down an expression for this operator. We will obtain a set of PDEs whose solutions are precisely the elements in its kernel, for which we can check that curves in the foliation which are nearby a fixed leaf correspond to solutions. By splitting the operator, and using automatic transversality \cite{Wen1}, we show that these are all possible solutions. This will imply that the index of the normal operator coincides with the dimension of its kernel, from which surjectivity follows. From the implicit function theorem, we also obtain regularity for Morse data chosen \emph{sufficiently close} to Morse-Bott data, which is enough for our purposes.

In order to do computations with linearized operators, we will choose a suitable connection on $W=\mathbb{R}\times M$. It will be given by the Levi--Civita connection of a suitable metric $g$ and hence symmetric. 

\paragraph*{Constructing a symmetric connection}

Given an almost complex structure $J$ which is compatible with a symplectic form $\omega$, we will denote by $g_{J,\omega}=\omega(\cdot,J\cdot)$ the associated Riemannian metric. 

\vspace{0.5cm}

Define, in the regions $\mathbb{R} \times Y\times \Sigma_\pm\backslash\mathcal{N}(\partial\Sigma_\pm)$, the metric 

$$g=\left(\begin{array}{cccc} 1  & 0 & 0 & 0 \\
0& g_{J_P^\pm,d\alpha_\pm}  &0 & 0 \\
0 & 0 & 1 & 0 \\
0 & 0 & 0 & g_{j_\pm,d\lambda_\pm} 
\end{array}\right),$$ where we are using the splitting $$T(\mathbb{R} \times Y\times \Sigma_\pm\backslash\mathcal{N}(\partial\Sigma_\pm))=\langle \partial_a \rangle \oplus \xi_\pm \oplus \langle R_\epsilon\rangle \oplus T\Sigma_\pm=\langle \partial_a \rangle \oplus TY \oplus T\Sigma_\pm$$

We extend it to $\mathbb{R}\times M_C^\pm$ by replacing $g_{j_\pm,d\lambda_\pm}$ by the identity in the basis $\{\partial_\rho,\partial_\theta\}$ in the above matrix. Along $\mathbb{R} \times Y\times I \times S^1$, set  
$$g=\left(\begin{array}{ccc} 1  & 0 & 0\\
0& g_{J_0,d\alpha} & 0\\
0 & 0 & 1  
\end{array}\right),$$ where we use the splitting $$T(\mathbb{R} \times Y\times I \times S^1)=\mathbb{R} \times T(Y \times I) \times \langle \partial_\theta\rangle$$

We set $\nabla=\nabla^g$ its Levi--Civita connection, which shall be the connection we will use to write down all our linearised Cauchy--Riemann operators.  

\subsubsection{Regularity for genus zero page-like curves in the Morse--Bott case}\label{reggenuszero}

In this section, we fix a genus zero curve $u:=u_{y,a}^-$ in the foliation, and denote by $u_{Y}:=u^{-1}(\mathbb{R}\times M_Y)$, $u_C:=u^{-1}(\mathbb{R}\times M_C^-)$, and $u_P:=u^{-1}(\mathbb{R}\times M_P^-)$. We will show that $\dim\ker \mathbf{D}_u^N=\mbox{ind}\;\mathbf{D}_u^N$, where $\mathbf{D}_u^N$ is the normal component of the linearized Cauchy-Riemann operator.

We will deal with the Morse--Bott case, where $H$ depends only on $r$. In this case, the operator we need to look at is given by 
$$
\mathbf{D}_u: W^{1,2,\delta_0}(u^*TW)\oplus  V_\Gamma \oplus X_\Gamma \rightarrow L^{2,\delta_0}(\overline{Hom}_{\mathbb{C}}((\dot{\Sigma}_-,j_-),(u^*TW,J))
$$
$$
\mathbf{D}_u\eta=\nabla \eta + J(u)\circ \nabla \eta \circ j_- + (\nabla_\eta J)\circ du \circ j_-
$$
Here, $\delta_0$ is a small weight making the operator Fredholm, $V_\Gamma$ is a $2k$-dimensional vector space of smooth sections asymptotic to constant linear combinations of $R_\epsilon$ and $\partial_a$, and $X_\Gamma$ is a $k(2n-1)$-dimensional vector space of smooth sections which are supported along $u_C\cup u_Y$ (a disjoint union of its cylindrical ends, which we denote $u_i$, $i=1,\dots, k$), and are constant equal to a vector in $T_yY$ along $u_Y$. We also have the operator
$$
\mathbf{D}_{(j_-,u)}\overline{\partial}_J: T_{j_-}\mathcal{T} \oplus W^{1,2,\delta_0}(u^*TW)\oplus V_\Gamma \oplus X_\Gamma \rightarrow L^{2,\delta_0}(\overline{Hom}_{\mathbb{C}}((\dot{\Sigma}_-,j_-),(u^*TW,J))
$$
$$
(Y,\eta)\mapsto J \circ du \circ Y + \mathbf{D}_u \eta,
$$ 
where $\mathcal{T}$ is a Teichm\"uller slice through $j_-$ (see e.g.\ \cite{Wen1} for a definition of this). The curve $u$ is said to be regular whenever $\mathbf{D}_{(j_-,u)}\overline{\partial}_J$ is surjective. By a result in \cite{Wen1}, this is equivalent to the surjectivity of its normal component. 
In this case, this operator is
$$
\mathbf{D}_u^N:W^{1,2,\delta_0}(N_u)\oplus X_\Gamma \rightarrow L^{2,\delta_0}(\overline{Hom}_{\mathbb{C}}(T\dot{\Sigma}_-,N_u))
$$
$$
\eta \mapsto \pi_N \mathbf{D}_u=\pi_N(\overline{\partial}_J\eta + (\nabla_\eta J)\circ du \circ j_-), 
$$
where $\pi_N$ is the orthogonal projection to the normal bundle $N_u$. Recall that the latter is any choice of $J$-invariant complement to the tangent space to $u$, which coincides with the contact structure at infinity. Riemann-Roch gives $\mbox{ind}\;\mathbf{D}_u^N=2n$. 

The operator $\mathbf{D}_u$ and $\mathbf{D}_u^	N$ are of Cauchy-Riemann type, which means in particular that they satisfy the Leibnitz rule 
\begin{equation}\label{Leib}
\mathbf{D}_u(f \eta)=f \mathbf{D}_u \eta+ \overline{\partial}f \otimes \eta
\end{equation}

\subparagraph*{Splitting over the paper}
      
We will think of the punctured surface $\dot{\Sigma}_-$ as being obtained abstractly from the surface with boundary $\Sigma_-\backslash\mathcal{N}(\partial\Sigma_-)$ by attaching cylindrical ends. Over the region $Y \times \Sigma_-\backslash\mathcal{N}(\partial\Sigma_-)$, we have a splitting $$u^*T(\mathbb{R}\times M)=T\Sigma_-\oplus T_{(a,y)}(\mathbb{R} \times Y)$$ 
Since $J$ preserves this splitting, this gives an identification $N_u=T_{(a,y)}(\mathbb{R}\times Y)=(\xi_-)_y\oplus\langle \partial_a, R_\epsilon(y)=R_-(y)/K \rangle$ of the normal bundle of $u$ along this region. Using that constant vectors in $N_u$ give holomorphic push-offs of $u$ in the foliation, we obtain    

$$\mathbf{D}_u=\left(
\begin{array}{cc}
\mathbf{D}_{id}\overline{\partial}_{j_-} & 0\\
0 & \mathbf{D}_u^N
\end{array}\right),$$

where the normal Cauchy--Riemann operator $\mathbf{D}_u^N$ splits as $\mathbf{D}_u^N=\bigoplus \overline{\partial}$.

\subparagraph*{Some technical assumptions} In order to be able to write down a manageable expression for $\mathbf{D}_u^N$ over the rest of the regions, we will assume, without loss of generality, that:
\begin{assumptions}\label{Assumptions}
\begin{enumerate}$\;$
\item[A.] $H$ has a unique critical point \emph{away} from a neighbourhood $(-\delta,\delta)\subset I$ where $J_0$ is non-cylindrical. Choose, say, $H^\prime(-2\delta)=0$.
\item[B.] Choose our coordinate $r$ so that the Liouville vector field $V$ coincides with $\pm \partial_r$ on the complement of $(-\delta,\delta)$.
\end{enumerate}
\end{assumptions}

These assumptions are only used in this section to show regularity, and do not affect other sections. Therefore we will, for simplicity, lift them in the rest of the sections.

\subparagraph*{Choosing a suitable normal bundle.} We now specify how we will extend our normal bundle $N_u$ to $u$ along its cylindrical ends. 

Since $J$ was chosen on $\xi_\epsilon\vert_{M_Y}$ so that the identification $(\xi_\epsilon,J)\rightarrow (T(Y\times I),J_0)$ is holomorphic, where $J_0$ may be \emph{any} $d\alpha$-compatible almost complex structure in $Y\times I$ which is cylindrical at the ends, we may identify the bundles. Observe that $J$ is always $\theta$-independent. And here we use assumptions A and B: we can choose $J_0$ so that it is \emph{cylindrical} in the complement of $(-\delta,\delta)$, so that $J$ is \emph{$r$-independent} ($r$ is the Liouville coordinate).

We then choose $N_u$ by the global expression
$$
N_u=\langle v,R_-(y)/K\rangle \oplus (\xi_-)_y=\langle v \rangle \oplus T_yY,
$$
where $v:=v_\rho= -J(R_-(y)/K)$ along the corner $M_C^-=Y\times (-\delta,\delta) \times S^1$, interpolating between $\partial_a=v_{-\delta}$ and $\partial_r=v_{\delta}$. Observe that $N_u$ is a trivial $J$-complex bundle.

\subparagraph*{Writing down the normal Cauchy-Riemann operator globally.} We now compute an asymptotic expression for $\mathbf{D}_u^N$. In the Morse-Bott case, one shows that $\mathbf{D}_u^N$ can be written asymptotically (i.e in the cylindrical coordinates $(s,t)$ along $u_{Y}$) as
$$
\mathbf{D}\eta(s,t)= \overline{\partial}\eta(s,t) + S_H(s,t)\eta(s,t), 
$$
where $S_H(s,t)=\nabla J$ is a symmetric matrix such that $S_H(s,\cdot)\rightarrow -2\pi \nabla_{-2\delta}^2H$ as $s\rightarrow +\infty$, uniformly in the second variable. 

By construction, $\nabla$ and $J$ are both independent of the coordinates $a,r,$ and $\theta$ along $\mathbb{R}\times M_Y$. Then
$$
S_H=\nu e_{11},
$$
where $e_{11}$ denotes the matrix with a $1$ at the $(1,1)$-entry, and zero everywhere else, and $\nu=-2\pi H^{\prime\prime}(-2\delta)>0$. As we traverse the smoothened corner $u_C$, we pick up a smooth path $(-\delta,\delta)\ni\rho \mapsto S_H(\rho)=s_H(\rho)e_{11}$, where $s_H(-\delta)=0$ and $s_H(\delta)=\nu$, so that 
$$
\mathbf{D}_u^N=\overline{\partial} + s_He_{11}
$$
\subparagraph*{Computing the kernel.} We write any section of $N_u$ as
$$
\eta=(\eta_0,\eta_\xi, \eta_\Gamma) \in W^{1,2,\delta_0}(\langle v, R_-(y)/K\rangle) \oplus W^{1,2,\delta_0}((\xi_-)_y) \oplus X_\Gamma, 
$$ 
with $\eta_0=(\eta_0^1,\eta_0^2) \in W^{1,2,\delta_0}(\langle v \rangle) \oplus  W^{1,2,\delta_0}(\langle R_-(y)/K\rangle)$. We denote 
$$
X_\Gamma^R:=\pi_R X_\Gamma, \; X_\Gamma^\xi:=\pi_\xi X_\Gamma,
$$
where 
$$\pi_R:T_yY \rightarrow \langle R_-(y)/K \rangle
$$
$$
\pi_\xi: T_yY \rightarrow (\xi_-)_y
$$ are the orthogonal projections with respect to the metric $g$. Then we write
$$
\eta_\Gamma=(\eta_\Gamma^R,\eta_\Gamma^\xi) \in X_\Gamma^R \oplus X_\Gamma^\xi,
$$ where $\eta_\Gamma^R=\pi_R \eta_\Gamma,$ $\eta_\Gamma^\xi=\pi_\xi \eta_\Gamma$. 

Then, $\eta \in \ker \mathbf{D}_u^N$ if and only if $\overline{\partial}\eta =-s_He_{11}\eta$, which in holomorphic coordinates $(s,t)$ is 
\begin{equation}\label{kerD}
\left\{\begin{array}{l} \partial_s\eta_0^1-\partial_t\left(\eta_0^2+\eta_\Gamma^R\right)=-s_H\eta_0^1\\
\partial_t\eta_0^1+\partial_s\left(\eta_0^2+ \eta_\Gamma^R\right)=0\\
\overline{\partial}\left(\eta_\xi + \eta_\Gamma^\xi\right)=0
\end{array}\right.
\end{equation}

Observe that the $\eta_\Gamma$ terms all disappear away from $u_C$. It is straightforward to check that all the nearby holomorphic curves in the foliation satisfy the above equations. We will show that these are indeed the unique solutions. 

We have shown that the operator $\mathbf{D}_u^N$ splits into a direct sum
$$
\mathbf{D}_u^N:=\mathbf{D}_u^\xi \oplus \mathbf{D}_u^R,
$$
where
$$
\mathbf{D}_u^\xi=\overline{\partial}: W^{1,2,\delta_0}((\xi_-)_y)\oplus X_\Gamma^\xi \rightarrow L^{2,\delta_0}(\overline{Hom}_{\mathbb{C}}(\dot{\Sigma}_-,(\xi_-)_y))
$$
and
$$
\mathbf{D}_u^R: W^{1,2,\delta_0}(\langle v, R_-(y)/K\rangle)\oplus X_\Gamma^R \rightarrow L^{2,\delta_0}(\overline{Hom}_{\mathbb{C}}(\dot{\Sigma}_-,\langle v, R_-(y)/K\rangle))
$$

We will show that both operators $\mathbf{D}_u^\xi$ and $\mathbf{D}_u^R$ are surjective, and this finishes the proof.

The first summand has kernel the $(2n-2)$-dimensional space of constant sections along $(\xi_-)_y$. Its index is $2n-2$, and it follows that it is surjective. The second summand satisfies
$$
\mathbf{D}_u^R=\overline{\partial}+s_H e_{11}
$$
In order to show that it is surjective, we use automatic transversality \cite{Wen1}. We need to check that $\mbox{ind }\mathbf{D}_u^R>c_1(\mathbf{D}_u^R)$, where $c_1(\mathbf{D}_u^R)$ denotes the \emph{adjusted first Chern number}, defined by 
$$
2c_1(\mathbf{D}_u^R)=\mbox{ind }\mathbf{D}_u^R-2+2g+\#\Gamma_{even},
$$
where $g$ is the genus, and $\Gamma_{even}$ is the set of punctures with even Conley-Zehnder index. 

The Conley-Zehnder index of $\mathbf{D}_u^R$ at each puncture is $1$, and therefore 
$$
\mbox{ind}\;\mathbf{D}_u^R=2-k+k=2 
$$
On the other hand, the adjusted first Chern number is
$$
2c_1(\mathbf{D}_u^R)=\mbox{ind}\;\mathbf{D}_u^R -2 + \#\Gamma_{even}=0
$$
This finishes the proof of regularity.

\subsection{From Morse--Bott to Morse}
\label{MBtoM}

In this section, we fix a Morse perturbation scheme. First, choose $H$ to be given by $H(y,r)=f(r)$ where $f$ is a sufficiently small and positive Morse function on $I$ and has a unique critical point at 0 with Morse index 1 (which yields the Morse--Bott situation). And then, choose a sufficiently small positive Morse function $g$ on $Y$ and extend it to a neighbourhood of $Y\times \{0\}$ to make a further perturbation (obtaining the Morse case). Therefore,
$$H(y,r)=f(r)+\gamma(r)g(y),$$  
where $\gamma:I\rightarrow [0,1]$ is a smooth bump function satisfying $\gamma=0$ in the region $\{|r|>2\delta\}$ and $\gamma=1$ on $\{|r|\leq\delta\}$. 

\begin{figure}[t] \centering
\includegraphics[width=0.25\linewidth]{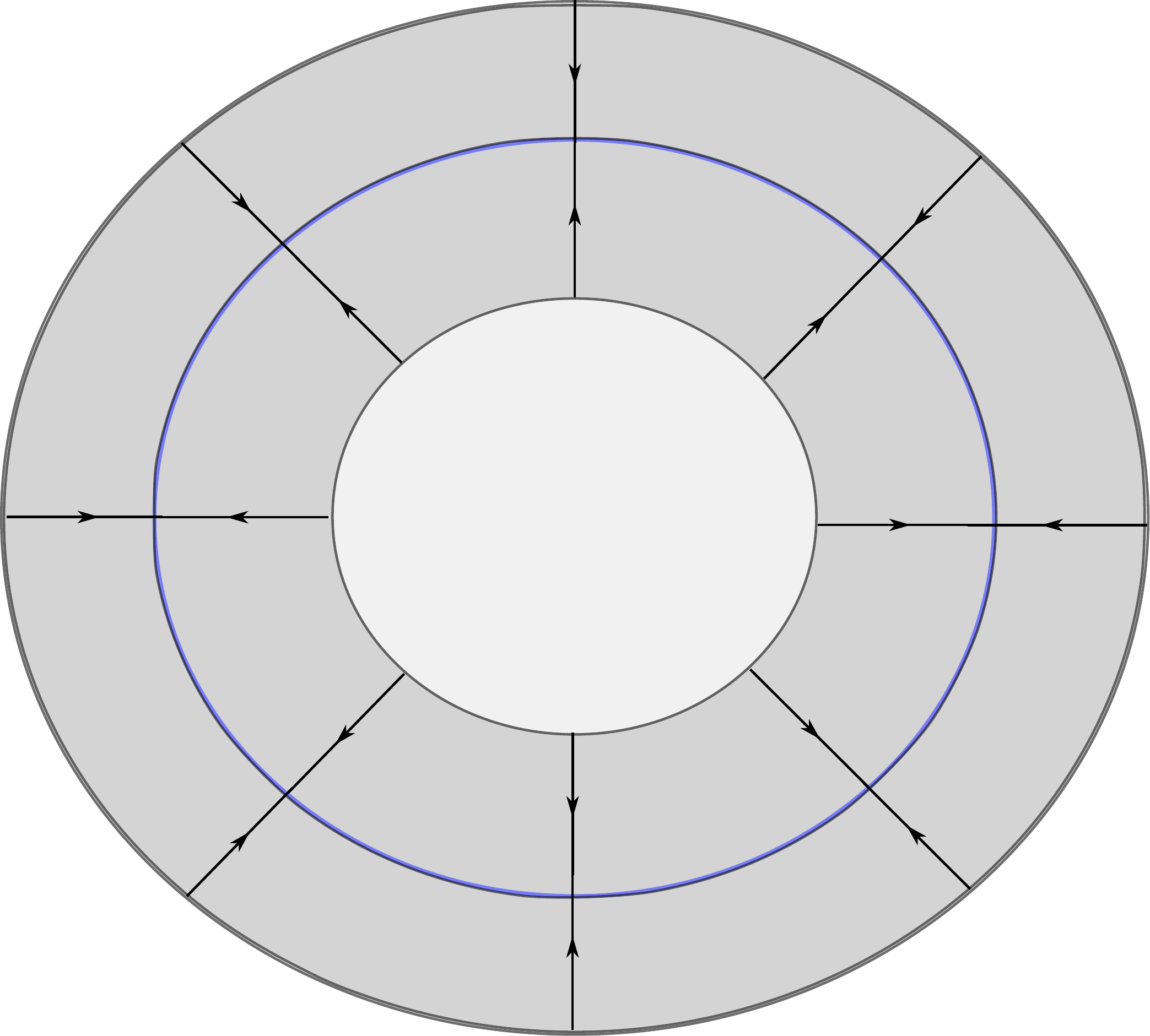}

\vspace{0.5cm}

\includegraphics[width=0.25\linewidth]{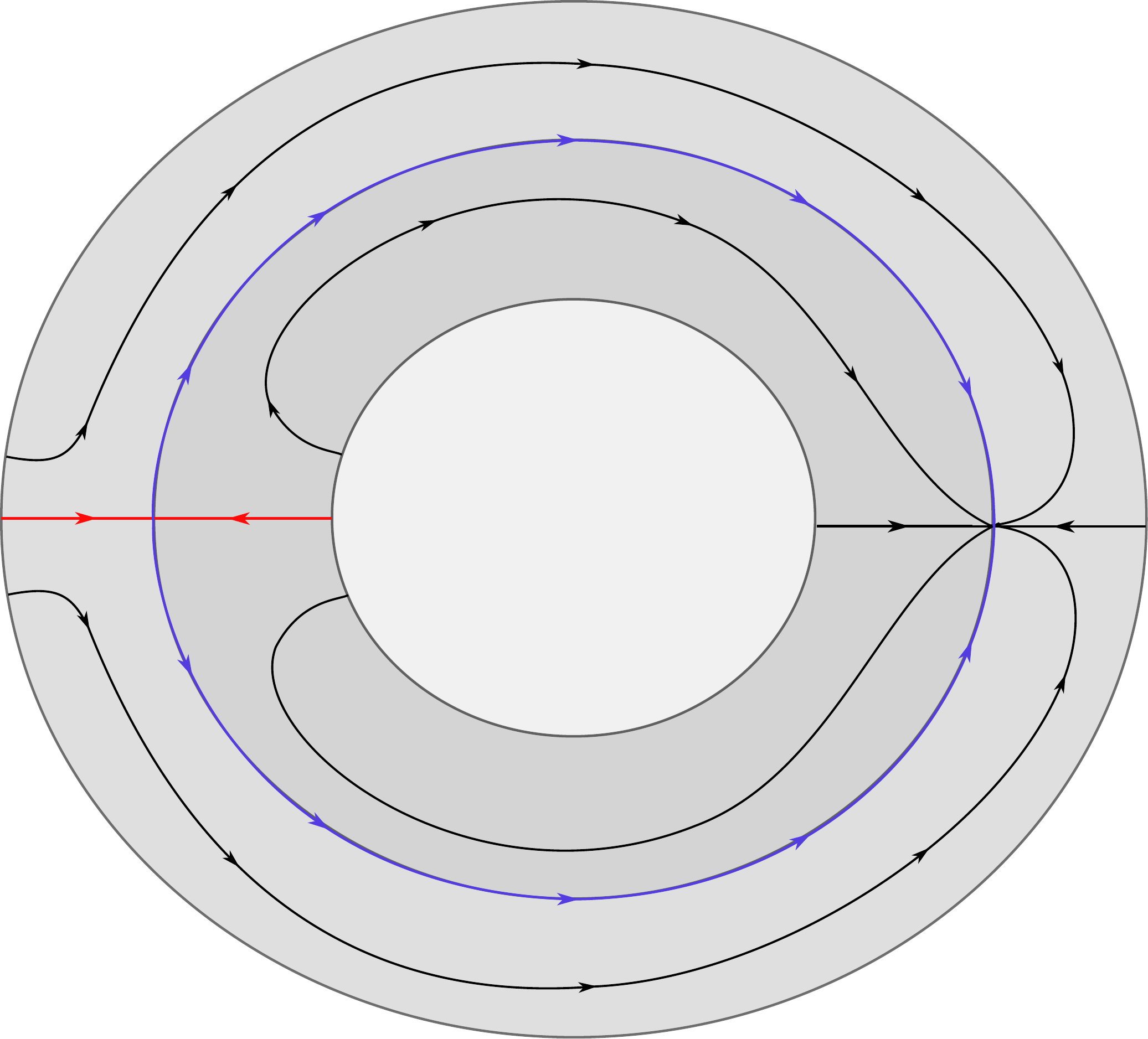}
\caption{\label{Morse} The Morse--Bott and the Morse scenarios, respectively, in the case where $Y=S^1$. The ``evil twins'' are shown in blue (see Section \ref{Torsion}).} 
\end{figure}

We view the Morse case as a deformation of the Morse--Bott one, via $H_t(y,r)=f(r)+t\gamma(r)g(y)$ for $t\in[0,1]$. We obtain a corresponding 1-parameter families of SHS's and compatible almost complex structures $J^t$. In the case where $g$ is chosen small, from the implicit function theorem we obtain:

\begin{thm}[Fredholm regularity in the nearby Morse case]
For Morse data sufficiently close to Morse-Bott data, all the genus zero curves in the finite energy foliation are Fredholm regular. 
\end{thm}

In order to simplify the torsion computation of Section \ref{Torsion}, we will choose $g$ to have a unique maximum $max$ and minimum $min$. Both scenarios are depicted in Figure \ref{Morse} in the case of $Y=S^1$.

\subsection{Uniqueness of the curves in the Morse/Morse--Bott case}
\label{uniquenessMB}

In this section, we prove that the family of curves we constructed above are the unique curves (up to reparametrization and multiple covers) that asymptote Reeb orbits of the family $\{\gamma_p: p \in \mbox{crit}(H)\}$, and with positive asymptotics in different components of $M_Y$. We do this in both the Morse/Morse--Bott situations. We assume that $H$ has a unique critical point in the interval direction, at $r=0$ (and perhaps other critical points contained in $Y\times \{0\}$, in the Morse case). In the Morse--Bott case, we denote by $\gamma_y:=\gamma_{(y,0)}=\{(y,0)\}\times S^1$ the simply covered Reeb orbit corresponding to $y \in Y$.

\begin{lemma}\label{posvsneg}
Assume either the Morse or Morse--Bott cases. Let $N\in \mathbb{N}^{>0}$, and let $u$ be a $J$-holomorphic curve with positive asympotics of the form $\gamma_p=\{p\}\times S^1$, for which the number of positive punctures is bounded above by $N$. Then, we can find sufficiently small $\epsilon>0$ (depending only on $N$), such that $$\#\Gamma^-(u)\leq \#\Gamma^+(u),$$ where $\Gamma^\pm(u)$ denotes the set of positive/negative punctures of $u$. Here, we count punctures with the covering multiplicity of their corresponding asymptotic.  
\end{lemma}
\begin{proof}
This follows easily from Remark \ref{remark11}, and by computing the $\Omega_\epsilon$-energy \cite{Mo2}. 
\end{proof}

\begin{thm}
\label{uniqueness}
Assume either the Morse or Morse--Bott scenario.

\vspace{0.5cm}

Let $u:\dot{\Sigma} \rightarrow \mathbb{R}\times M$ be a (not necessarily regular) $J$-holomorphic curve defined on some punctured Riemann surface $\dot{\Sigma}$ which is asymptotically cylindrical, and asymptotes simply covered Reeb orbits of the form $\gamma_p=\{p\}\times S^1$ at its positive ends. Assume that any two of the positive ends of $u$ lie in distinct components of $M_Y$. 

\vspace{0.5cm}

Then, for sufficiently small and uniform $\epsilon>0$, we have that, if $u$ is not a trivial cylinder over one of the $\gamma_p$'s, then $u$ is a curve of the form $u_{y,a}^\pm$ for some $y \in Y$, or a flow-line cylinder $u_\gamma$ (in the Morse scenario).
\end{thm}

\begin{proof}

We will consider two cases: either $u$ is completely contained in the region $\mathbb{R}\times M_Y$ (case A), or it is not (case B). 

\vspace{0.5cm}

\textbf{Case A.} This case is easily dealt with in the Morse--Bott scenario.  By assumption, we have that $u$ has a unique positive end. Since the $\gamma_p$'s are not contractible/nullhomologous inside $M_Y$, Lemma \ref{posvsneg} implies that $u$ is has one positive and one negative end, both simply covered, corresponding to Reeb orbits $\gamma_{p_\pm}$. But then the $\Omega_\epsilon$-energy of $u$ vanishes, and $u$ is necessarily is a cover of a trivial cylinder.

\vspace{0.5cm} 

In the Morse case, we show that $u$ is a flow line cylinder. Again, Lemma \ref{posvsneg} implies that $u$ is has one positive and one negative end, both simply covered, corresponding to Reeb orbits $\gamma_{p_\pm}$. Observe that, a priori, $u$ is not even necessarily a cylinder, since it may have positive genus. 

In the degenerate case when $H\equiv 0$, we the projection $\pi_Y:\mathbb{R}\times M_Y \rightarrow Y \times I$ is holomorphic. Then, if $u$ is a holomorphic curve for this data, then so is $v:=\pi_Y \circ u$. Since the asymptotics of $u$ are covers of the $S^1$-fibers of $\pi_Y$, the map $v$ extends to a holomorphic curve in the closed surface $\Sigma$. But $Y\times I$ is exact, so that $v$ has to be constant by Stokes' theorem. This means that $u$ is necessarily a multiple cover of a trivial cylinder.

We then see that the space of \emph{stable holomorphic cascades} \cite{Bo} in $\mathbb{R}\times Y \times I \times S^1$, the objects one obtains as limits of honest curves when one turns off the function $H$, consists of finite collections of flow-line segments and covers of trivial cylinders. If we take $H_t=tH$, and we assume we have a sequence $\{u_n\}$ of $J_{t_n}$-holomorphic maps with $t_n\rightarrow 0$ (where $J^{t}$ is the almost complex structure corresponding to $H_{t}$), with one positive and negative simply covered orbits corresponding to critical points $p_\pm$, then we obtain a \emph{stable} holomorphic cascade $\mathbf{u}^H_\infty$ as a limiting object. Since the positive end of $u_n$ is simply covered, Lemma \ref{posvsneg} applied to $H=0$ implies that \emph{every} Reeb orbit appearing in $\mathbf{u}^H_\infty$ is simply covered, and therefore every of its holomorphic map components cannot be multiply covered. These can then only be trivial cylinders, but stability of the cascade means that it does not have trivial cylinder components. We conclude that the space of holomorphic cascades which glue to curves as in our hypothesis consists solely of flow-lines, which are regular by the Morse--Smale condition, and come in a $(ind_{p_+}(H)-ind_{p_-}(H)-1)$-dimensional family. The gluing results in \cite{Bo} give a 1-1 correspondence between index $d$ Morse families of curves for the non-degenerate perturbation and index $d$ regular holomorphic cascades, which in our situation is exactly what proves that our curve $u$ is a flow-line cylinder. 

\vspace{0.5cm}

\textbf{Case B.} We first assume the Morse--Bott case, and deal with the Morse case via the gluing results in \cite{Bo}. The approach in this situation is to estimate the $\Omega_\epsilon$-energy of $u$, and to use a suitable branched cover argument. The details are as follows.  

\vspace{0.5cm}

Assume the Morse--Bott case. Since every positive puncture of $u$ corresponds to a critical point, Remark \ref{remark11} implies that so does every negative one. Let us denote by $\Gamma^\pm$ the set of positive/negative punctures of $u$, and for $z \in \Gamma^\pm$, let $\gamma^{\kappa_z}_{p_z}$ be the Reeb orbit corresponding to $z$, where $p_z\in \mbox{crit}(H)$, $\kappa_z\geq 1$, and $\kappa_z=1$ for $z \in \Gamma^+$. By assumption, we have that $\#\Gamma^+\leq k$, the number of components of $M_Y$. 

The $\Omega_\epsilon$-energy then has the following upper bound:

\begin{equation}
\label{firstineq}
\begin{split}
\mathbf{E}(u)=&\int_{\dot{\Sigma}}u^*\Omega_\epsilon\\
=&2\pi e^{\epsilon H(0)}\left(\#\Gamma^+-\sum_{z\in\Gamma^-}\kappa_z\right)\\
=&2\pi e^{\epsilon H(0)}(\#\Gamma^+-\#\Gamma^-)\\
\leq&2\pi\#\Gamma^+||e^{\epsilon H}||_{C^0}\\
\leq& 2\pi k ||e^{\epsilon H}||_{C^0}
\end{split}
\end{equation} 

By construction, we have that over the region $\Sigma_\pm\backslash\mathcal{N}(\partial \Sigma_\pm)$, the almost complex structure $J$ splits. This implies that the projection $$\pi_P^\pm:\mathbb{R}\times Y\times \Sigma_\pm\backslash\mathcal{N}(\partial \Sigma_\pm) \rightarrow\Sigma_\pm\backslash\mathcal{N}(\partial \Sigma_\pm)$$ is holomorphic.
Moreover, this is still true if we extend this region by adding a small collar $\{t=\mp \rho \in (-\delta, -\delta+\delta/9]\}$, as one can check. Denote by $B_\pm^\delta:=\Sigma_\pm\backslash\mathcal{N}(\partial \Sigma_\pm)\cup\{t \in(-\delta,-\delta+\delta/9]\}$.

For each $w \in B_\pm^\delta$, the hypersurface $E_w:=\mathbb{R}\times Y \times \{w\}$ is $J$-holomorphic, and $\pi_P^\pm$ has $E_w$ as fiber over $w$. Since the asymptotics of $u$ are away from $B_\pm^\delta$, the intersection of $u$ with any of the $E_w$ is necessarily a finite set of points, since they are restricted to lie in a compact part of the domain of $u$.

Assuming WLOG that $u$ indeed has a non-empty portion lying over the ``plus'' region $\Sigma_+\backslash{\mathcal{N}(\partial \Sigma_+)}$, by positivity of intersections we have that $u$ necessarily intersects every $E_w$ for $w\in B_\pm^\delta$. By Sard's theorem, we may then find a $t_0\geq -\delta$, so that $\pi_P^+(u)$ is transverse to the circle $\{t_0\} \times S^1 \subseteq B_+^\delta$ (over all $k$ components of the collar). If we denote $B_+^{\delta;t_0}:= B_+^\delta\backslash\{t \in (t_0,-\delta+\delta/9]\}$, we have that $S_+:=(\pi_P^+\circ u)^{-1}(B_+^{\delta;t_0})\subseteq \dot{\Sigma}$ is a surface with boundary $\partial S_+=(\pi_P^+\circ u)^{-1}(\partial B_+^{\delta;t_0})$. The upshot of the discussion above is that the map
$$
F_u:=(\pi_P^+\circ u)\vert_{S_+}:S_+ \rightarrow B_+^{\delta;t_0}
$$ 
is a \emph{holomorphic branched cover}, having as degree the (positive) algebraic intersection number of $u$ with any of the $E_w$ (which is independent of $w \in B_+^{\delta;t_0}$). Call this degree $\deg^+(u):=\deg(F_u)$. We wish to show that $\deg^+(u)=1$, and so this map will be actually a biholomorphism.

\vspace{0.5cm}

Let us write $\partial S_+=\bigcup_{i=1}^l C_i$, where $C_i$ is a simple closed curve whose image under $F_u$ wraps around one of the circles $\{t_0\}\times S^1$, with winding number $n_i$. By holomorphicity of $F_u$, we have $n_i>0$. Observe that necessarily one has that $l\geq k$, since $u$ intersects every $E_w$ at least once, and in particular for every $w$ on the $k$ circles $\{t_0\}\times S^1$. 

By counting preimages under this projection of a point in each the $k$ circles $\{t_0\}\times S^1$, we obtain that

\begin{equation}\label{degreeplus}
\int_{\partial S_+} u^*d\theta = 2\pi\sum_{i=1}^l n_i=2\pi k \deg^+(u)
\end{equation}      

Using expression $\Omega_\epsilon=Kd\alpha_\pm+d\lambda_\pm$, equation (\ref{degreeplus}), the fact that $t_0\geq-\delta$, $u^*d\alpha_+\geq 0$, and Stokes' theorem, we have the following energy estimate:

\begin{equation}
\label{ineqs}
\begin{split}
\mathbf{E}(u)\geq& \int_{S_+}u^*\Omega_\epsilon\\
=&\int_{S_+}u^*\left(K d\alpha_++ d\lambda_+\right)\\
\geq&\int_{S_+}u^*d\lambda_+\\
=&\int_{\partial S_+}u^*(e^{t}d\theta)\\
=&2\pi k e^{t_0}\deg^+(u)\\
\geq & 2\pi k e^{-\delta}\deg^+(u)
\end{split}
\end{equation}

If we combine this with the inequality (\ref{firstineq}), we obtain 
\begin{equation}
\label{sss}
2\pi k ||e^{\epsilon H}||_{C^0} \geq 2\pi e^{\epsilon H(0)}\left(\#\Gamma^+-\#\Gamma^-\right)\geq 2\pi k e^{-\delta} \deg^+(u) \geq 2\pi k e^{-\delta}
\end{equation}
Now, since we have the freedom to choose $-\delta$ and $||\epsilon H||_{C^0}$ as close to zero as wished, one can easily see that $$\#\Gamma^+=k,\;\#\Gamma^-=0,\deg^+(u)=1$$ 

This proves that $u$ has no negative ends and precisely $k$ positive ends,  and that $F_u$ gives a biholomorphism between $S_+$ and $B_+^{\delta;t_0}$. It also follows from our energy estimates that $$\int_{S_+}u^*d\alpha_+=0$$ Since the integrand is non-negative, we get that $u^*d\alpha_+=0$, so that $u(S_+)$ lies entirely in the almost complex 4-manifold 
$$
U:=\mathbb{R}\times \gamma \times\Sigma_+\subseteq \mathbb{R}\times Y \times \Sigma,
$$
where $\gamma$ is a (not necessarily closed) $R_+$-Reeb orbit. 

Choose now a point $y$ in the projection of $u(S_+)$ to $\gamma$, in this region. We shall prove that the projection of $u(S_+)$ to $\gamma$ consists of just the point $y$, using the Morse--Bott assumption. 

Using the Reeb flow along $\gamma$ we have a local coordinate $y=y(s)$, such that $y(0)=y$. Assume by contradiction that there is some $s\neq 0$ such that $y(s)$ belongs to the projection of $u(S_+)$ to $\gamma$, and consider the family of curves $\{u_{y(s),a}^+:a\in \mathbb{R}\}$. By choosing a suitable $a$, we obtain an intersection of $u_{y(s),a}^+$ and $u$, and, by positivity of intersections, we may assume that $y(s)$ is such that $\gamma_{y(s)}$ is not an asymptotic orbit of $u$. Then the total intersection of these curves in $U$ is at least 1.

On the other hand, the set $\mathbb{R} \times (u_M \cap (u_{y(s),a}^+)_M)$, where $u_M$ is the projection to $M$ of $u$, must be bounded in the $\mathbb{R}$-direction, since the corresponding Reeb orbits are bounded away from each other. This means that the standard intersection pairing is homotopy invariant (there are no contributions coming from infinity), but a priori the intersection points might move along with a given homotopy. If we choose to homotope by translating in the $\mathbb{R}$-direction, this can only happen if the projection of both curves to $M$ intersect the asymptotics of each other. Since the projection of the curve $u_{y(s^\prime),a^\prime}$ to $M$ does not intersect the asymptotics of $u$, we can homotope the intersections away, a contradiction. This proves the claim that $u$ is constant in $Y$.

We have obtained that the portion of $u$ which lies in the over $\Sigma_+\backslash \mathcal{N}(\partial \Sigma_+)$ is actually contained in the 3-manifold $M_y:=\mathbb{R}\times \{y\} \times \Sigma_+\backslash \mathcal{N}(\partial \Sigma_+)$, as is the corresponding portion of the curve $u_{y,a}^+$. But for dimensional reasons, this manifold has a unique 2-dimensional $J$-invariant distribution, given by $TM_y\cap JTM_y$. Therefore the tangent space to $u$ must coincide with the tangent space of some $u_{y,a}^+$ at every point in this region. The unique continuation theorem finally yields $u=u_{y,a}^+$. 

\vspace{0.5cm}

This proves the theorem in the Morse--Bott situation. The proof in the Morse case follows from uniqueness in the Morse--Bott one, and the gluing results in \cite{Bo}.

\end{proof}

\subsection{From the SHS to a sufficiently non-degenerate contact structure}
\label{SHStoND}

For computations in SFT we need non-degenerate Reeb orbits and contact data. Therefore, we need to perturb the SHS $\mathcal{H}_\epsilon=(\Lambda_\epsilon,\Omega_\epsilon)$ to a nearby contact structure.

\subparagraph*{Perturbation to contact data} Recall that we have defined an exact symplectic form $\omega_\sigma$ on the double completion, and we denote by $V_\sigma$ its associated Liouville vector field. We also have the ``vertical'' Liouville vector field $X$ associated to $d\lambda$, defined by expression (\ref{X}), and the stabilizing vector field $Z$ for the SHS $\mathcal{H}_\epsilon$, defined by (\ref{Z}).

For $s \in [0,1]$, define 
$$Z_s=\left\{\begin{array}{ll} V_\sigma+((1-s)\beta(t)+s)X, & \mbox{ in the region } E^{\infty,\infty}(t) \mbox{ where } t \mbox{ is defined}\\ 
\frac{\sigma}{\sigma+\sigma^\prime}\partial_r+sX_+,& \mbox{ in } E^{\infty,\infty}(t)^c \cap \{r > 1-\delta\}\\
-\frac{\sigma}{\sigma-\sigma^\prime}\partial_r+sX_-,& \mbox{ in } E^{\infty,\infty}(t)^c \cap \{r < -1+\delta\}
\end{array}\right. $$ 

We have that $Z_1=X_{\sigma}$ is the Liouville vector field associated to $\omega_{\sigma}$ and $Z_0=Z$. This yields a family of SHS's given by $$\mathcal{H}_{\epsilon,s}=(\Lambda_{\epsilon,s},\Omega_{\epsilon})_{s\in[0,1]}=(i_{Z_s}\omega_{\sigma}\vert_{TM^\epsilon},\omega_{\sigma}\vert_{TM^\epsilon})_{s\in[0,1]}$$ 

One can see that $\Lambda_{\epsilon,s}$ is a contact form for $s>0$. By Gray's stability, as long as $\epsilon$ and $s$ are positive and sufficiently close to zero, the isotopy class of $\xi_{\epsilon,s}$ is independent on parameters.

\subparagraph*{Holomorphic curves for the contact data} Since we have shown that the genus zero holomorphic curves are regular for the SHS data, the implicit function theorem implies that they will survive a small perturbation to contact data, and will still be regular. 

\subsection{Proof of Theorem \ref{thm5d}}
\label{Torsion}

Once all the technical tools are in place, we will prove Theorem \ref{thm5d}. We fix the parameters $\epsilon$ and $s$, so that we work in the SFT algebra $\mathcal{A}(\Lambda_{\epsilon,s})[[\hbar]]$ whose homology is $H_*^{SFT}(M,\xi_{\epsilon,s})$ (which is independent on parameters). We take coefficients in $R_{\Omega}=\mathbb{R}[H_2(M;\mathbb{R})/\ker \Omega]$, where $\Omega \in \Omega^2(M)$ defines an element in the annihilator $\mathcal{O}:=$Ann$(\bigoplus^k H_1(Y;\mathbb{R})\otimes H_1(S^1;\mathbb{R}))$. Here, we view $\bigsqcup^k Y\times S^1=\bigsqcup^k Y\times S^1\times \{0\}$ as sitting in $M_Y\subseteq M$. We will show that $(M,\xi_{\epsilon,s})$ has $\Omega$-twisted $(k-1)$-torsion for every $[\Omega] \in \mathcal{O}$, so that in particular it has $(k-1)$-torsion for \emph{untwisted} version of the SFT algebra.

Let us recall that for SFT to be defined, we need to introduce an \emph{abstract perturbation} of the Cauchy--Riemann equation. We shall be doing our computation \emph{prior} to introducing this perturbation, and prove by the end the section that this is a reasonable thing to do. See the end of this section for more details.

\subparagraph*{Computation of algebraic torsion} The index formula (\ref{index}) implies that curves $u_{y,a}^-$ which asymptote index $2n$ critical points (maxima) at $k-1$ of the $k$ positive ends, and one index $1$ critical point at the remaining one, have index $1$. 

Given $e_1,\dots,e_{k-1}$ maxima and $h$ an index $1$ critical point, denote by $\gamma_{e_i}$ and $\gamma_h$ the corresponding non-degenerate Reeb orbits. Consider the moduli space $$\mathcal{M}=\mathcal{M}_{0}(W,J_{\epsilon,s};(\gamma_{e_1},\dots, \gamma_{e_{k-1}},\gamma_h),\emptyset)$$ of $\mathbb{R}$-translation classes of $k$-punctured genus zero $J_{\epsilon,s}$-holomorphic curves in $W=\mathbb{R}\times M$, which have no negative asymptotics, and have $\gamma_{e_1},\dots,\gamma_{e_{k-1}},\gamma_h$ as positive asymptotics. The uniqueness Theorem \ref{uniqueness} implies that every element in $\mathcal{M}$ is a genus zero curve in our foliation. Since every such curve is regular, after a choice of coherent orientations as in \cite{BoMon}, this moduli space is an oriented zero dimensional manifold, which can therefore be counted with appropiate signs. Now, a choice of coherent orientations for the moduli space of Morse flow lines induces a coherent orientation for the moduli space containing the curves $u^a_\gamma$, and we fix such a choice here onwards. We choose our function $$H:Y\times I \rightarrow \mathbb{R}^{\geq 0}$$ as made explicit in Section \ref{MBtoM}. In particular, there are no index zero critical points, and the only index 1 critical point is given by $(min,0)$, where $min$ is the unique minimum of $g$. We shall denote by $\mathcal{M}(H;p_-,p_+)$ the moduli space of \emph{positive} unparametrized flow lines $\gamma$ connecting $p_-$ to $p_+$. Assuming the Morse--Smale condition for $H$, we have then that the zero dimensional moduli spaces correspond to critical points satisfying $ind_{p_+}(H)=ind_{p_-}(H)+1$. 

\vspace{0.5cm}

We then fix  $e_1,\dots,e_{k-1},h=(min,0)$ as above, and we let $q_{e_i}$ and $q_h$ be the generators in SFT corresponding to the Reeb orbits associated to these critical points. Let $$Q=q_{e_1}\dots q_{e_{k-1}}q_h,$$ which is an element of $\mathcal{A}(\Lambda_{\epsilon,s})[[\hbar]]$. In order to compute its differential, we need to count all of the rigid holomorphic curves which asymptote at its positive ends any of the Reeb orbits appearing in $Q$.

\vspace{0.5cm}

\textbf{Claim}. The holomorphic curves contributing to the differential of $Q$ are of either the two following types: 

\begin{itemize}
\item A holomorphic sphere $u_{y,a}^- \in \mathcal{M}$, which is in fact unique.

\item A holomorphic cylinder $u_\gamma$ inside $r=0$, connecting an index $2n-1$ critical point to one of the maxima $e_i$.
\end{itemize}

Indeed, using Theorem \ref{uniqueness}, it follows that only the somewhere injective curves in our foliation are involved in the computation of this differential (see Proposition \ref{buildingsuni} below). Using that there are no index zero critical points, we see that there are only two possible ways to approach the critical point $(min,0)$, entering through the two different boundary components of $Y \times I$, and that there is a unique $\mathbb{R}$-class of the form $u_y^-$ which has $\gamma_h$ as a positive asymptotic (see Figure \ref{Morse}). Moreover, by observing that the generic behaviour is hitting a maxima, by a generic and small perturbation of the Morse function $H$ along different components of the spine, we can arrange that this curve actually defines an element of $\mathcal{M}$. The uniqueness Theorem \ref{uniqueness} above shows that in fact this is the only element in $\mathcal{M}$. Finally, ruling out the curves coming from the positive side (which have the wrong index and hence are not counted by SFT), we are left only with the curves listed in our claim.

\vspace{0.5cm}

In order to count the curves of type $u_\gamma$, we observe that positive flow lines going from an index $2n-1$ critical point $p$ to an index $2n$ critical point $q$ come in ``evil twins'' pairs $\gamma \leftrightarrow \overline{\gamma}$: by definition of the Morse index, we have only one positive eigendirection for the Hessian of the Morse function at $p$, and the flow lines approach this point on either side of this direction. Since $H_{2n}(Y\times I)=H_{2n}(Y)=0$ ($Y$ being a $(2n-1)$-manifold), we only have one generator of the top Morse homology chain group, which is necessarily closed under the Morse differential.  Therefore, after choosing a coherent orientation of the moduli spaces of curves, the evil twins cancel each other out, and hence $\# \mathcal{M}(H;p,q)=0$. 

\vspace{0.5cm}

Fix $\Omega$ a closed $2$-form so that $[\Omega] \in \mathcal{O}$. For $d \in H_2(M;\mathbb{R})$, we denote by $\overline{d} \in H_2(M;\mathbb{R})/\ker \Omega$ the class it induces, by $u_0$ the unique $\mathbb{R}$-translation class in $\mathcal{M}$, and, for a rigid holomorphic curve $v$, we denote by $\epsilon(v)$ the sign of $v$ assigned by our choice of coherent orientations. In particular, we know that $\epsilon(u_\gamma)=-\epsilon(u_{\overline{\gamma}})$. 

Observe that, for the $i$-th component of $M_Y$, Reeb orbits corresponding to critical points all define the same homology class $[S^1]_i \in H_1(M)$. We take as canonical representative of this class the $1$-cycle given by the Reeb orbit $\gamma_{e_i}$ over the unique maxima $e_i$. For every index $2n-1$ critical point $p$ lying in the $i$-th component, fix $\gamma$ an index $1$ flow line joining $p$ to the maxima $e_i$. Choose the spanning surface of $\gamma_p$ to be $F_p:= -\gamma \times S^1$, satisfying $\partial F_p = \gamma_p-\gamma_{e_i}$. Then, for this choice of spanning surfaces, the homology class associated to the holomorphic cylinder $u_\gamma$ is $[u_\gamma]=[F_p\cup u_\gamma]$. One thinks of $F_p$ as being attached to $u_\gamma$ at the negative end, corresponding to $\gamma_p$. Let $T_{\gamma,\overline{\gamma}}$ denote the $2$-torus $\overline{\gamma\cup \overline{\gamma}}\times S^1$. Observe that $$[u_\gamma]-[u_{\overline{\gamma}}]=[T_{\gamma,\overline{\gamma}}] \in H_1(Y)\otimes H_1(S^1)\subseteq \ker \Omega$$ 
Therefore, we have $\overline{[u_\gamma]}=\overline{[u_{\overline{\gamma}}]}$.

According to Proposition \ref{buildingsuni} below, the image of $Q$ under the differential $\mathbf{D}_{\epsilon,s}:\mathcal{A}(\Lambda_{\epsilon,s})[[\hbar]] \rightarrow \mathcal{A}(\Lambda_{\epsilon,s})[[\hbar]]$ is then given by 

\begin{equation}
\begin{split}
\mathbf{D}_{\epsilon,s}(Q)=&\epsilon(u_0)z^{\overline{[u_0]}}\hbar^{k-1} + \sum_{\substack{i=1,\dots,k-1\\ind_p(H)=2n-1}}\sum_{\gamma \in \mathcal{M}(H;p,e_i)}\epsilon(u_\gamma)z^{\overline{[u_\gamma]}}q_{p}\frac{\partial Q}{\partial q_{e_i}}\\
=&\epsilon(u_0)z^{\overline{[u_0]}}\hbar^{k-1} + \frac{1}{2}\sum_{\substack{i=1,\dots,k-1\\ind_p(H)=2n-1}}\sum_{\gamma \in \mathcal{M}(H;p,e_i)}(\epsilon(u_\gamma)+\epsilon(u_{\overline{\gamma}}))z^{\overline{[u_\gamma}]} q_{p}\frac{\partial Q}{\partial q_{e_i}}\\
=& \epsilon(u_0)z^{\overline{[u_0]}}\hbar^{k-1},
\end{split}
\end{equation}
which proves that our model has $(k-1)$-torsion with coefficients in $R_{\Omega}$. We have used that all orbits are simply covered, so no combinatorial factors appear, and we need not worry about asymptotic markers.

\subparagraph*{Why the computation works} We now justify the computation above. Let us recall first the fact that the abstract polyfold machinery for SFT requires the introduction of an abstract perturbation to the Cauchy--Riemann equation making every holomorphic curve of positive index regular. The basic facts about this perturbation scheme, which comes from the polyfold theory of Hofer--Wysocki--Zehnder, are that 

\begin{itemize}
\item Every Fredholm regular index 1 holomorphic curve gives rise to a unique solution to the perturbed problem, if the perturbation is sufficiently small.
\item If solutions to the perturbed problem with given asymptotic behaviour exist for small perturbations, then as the perturbation is switched off they give rise to a subsequence of curves which converge to a holomorphic building with the same asymptotic behaviour. 
\end{itemize}

Therefore, the index 1 curves in our foliation survive and are counted, but we need to make sure there are no extra curves which need to be taken into the count. In what follows, $J$ will denote the original $J_{\epsilon,0}$ we constructed in the Morse case. 

\begin{prop}\label{buildingsuni} The space of connected $J$-holomorphic stable buildings of index 1 which may become curves contributing to the differential of $Q=q_{e_1}\dots q_{e_{k-1}}q_h$ after introducing an abstract perturbation, or after perturbing the SHS to a sufficiently nearby contact structure, have only one level, no nodes, and are somewhere injective and regular: they actually consist exactly of either an index 1 curve $u_\gamma$ for some positive flow line $\gamma$, or a curve $u^-_{y,a}$ for some $y \in Y$, $a \in \mathbb{R}$.
\end{prop}

\begin{proof} It follows by inductively applying Theorem \ref{uniqueness} \cite{Mo2}.
\end{proof}

\section{Giroux torsion implies algebraic 1-torsion in higher dimensions}
\label{GirouxTorsion1-tors}

In this section, we address a conjecture in the paper \cite{MNW} (Conjecture 4.14). 

\begin{defn}[Giroux]\label{idLD} Let $\Sigma$ be a compact $2n$-manifold with boundary, $\omega$ a symplectic form on the interior $\mathring{\Sigma}$ of $\Sigma$, and $\xi_0$ a contact structure on $\partial \Sigma$. The triple $(\Sigma,\omega,\xi_0)$ is an \emph{ideal Liouville domain} if there exists an auxiliary 1-form $\beta$ on $\mathring{\Sigma}$ such that
\begin{itemize}
\item $d\beta = \omega$ on $\mathring{\Sigma}$.
\item For any smooth function $f:\Sigma \rightarrow [0,\infty)$ with regular level set $\partial\Sigma=f^{-1}(0)$, the 1-form $\beta_0=f\beta$ extends smoothly to $\partial \Sigma$ such that its restriction to $\partial \Sigma$ is a contact form for $\xi_0$.
\end{itemize} 
The 1-form $\beta$ is called a \emph{Liouville form} for $(\Sigma,\omega,\xi_0)$.
\end{defn}

In \cite{MNW}-terminology, we say that an oriented hypersurface $H$ in a contact manifold $(M,\xi)$ is a \emph{$\xi$-round hypersurface} modeled on some closed contact manifold $(Y, \xi_0)$ if it is transverse to $\xi$ and admits an orientation preserving identification with $S^1\times Y$ such that $\xi \cap TH = TS^1 \oplus \xi_0$. 

Given an ideal Liouville domain $(\Sigma,\omega,\xi_0)$, the \emph{Giroux domain} associated to it is the contact manifold $\Sigma \times S^1$ endowed with the contact structure $\xi=\ker(f(\beta + d\theta))$, where $f$ and $\beta$ are as before, and $\theta$ is the $S^1$-coordinate. Away from $V(f)=\partial \Sigma$, the vanishing locus of $f$, this contact structure coincides with the \emph{contactization} $\ker(\beta + d\theta)$. Over $V(f)$ it is just given by $\xi_0 \oplus TS^1$, so that $V(f)\times S^1$ is a $\xi_0$-round hypersurface modeled on $(\partial \Sigma,\xi_0)$.

One may find a collar neighbourhood of the form $[0,1)\times H\subseteq \Sigma \times S^1$, on which $\xi$ is given by the kernel of a contact form $\beta_0+sd\theta$, where $s$ is the coordinate on the interval, where $H \subseteq \partial M$ corresponds to $s=0$, $\theta$ the coordinate in $S^1$, and $\beta_0$ a contact form for $\xi_0$ \cite[Lemma 4.1]{MNW}. Using these collar neighbourhoods, one has a well-defined notion of gluing of two Giroux domains along boundary components modeled on isomorphic contact manifolds (see Section \ref{GirouxSOBDs} below). 

We also have a blow-down operation for round hypersurfaces lying in the boundary. If $H$ is a $\xi$-round boundary component of $(M, \xi)$, with orientation opposite
the boundary orientation, consider the collar neighborhood $[0,1)\times H$ as before. Let $\mathbb{D}$ be the disk of radius $\sqrt{\epsilon}$ in $\mathbb{R}^2$. The map $\Psi: (re^{i\theta},y) \mapsto (r^2,\theta, y)$ is a diffeomorphism from $(\mathbb{D} \backslash \{0\}) \times Y$ to $(0, \epsilon) \times S^1 \times Y$ which pulls back $\beta_0 + s dt$ to the contact
form $\beta_0 + r^2d\theta$. Thus we can glue $\mathbb{D} \times Y$ to $M \backslash H$ to get a new contact manifold in which $H$ has been replaced by $Y$, and the $S^1$-component of $H$ has been capped off.

Given a contact embedding of the interior of a Giroux domain $G_\Sigma
:=\Sigma\times S^1$ inside a contact manifold $(M,\xi)$, we shall denote by $\mathcal{O}(G_\Sigma)\subseteq H^2(M;\mathbb{R})$ the annihilator of $H_1(\Sigma;\mathbb{R})\otimes H_1(S^1;\mathbb{R})$, when the latter is viewed as a subspace of $H_2(M;\mathbb{R})$. If $N\subseteq (M,\xi)$ is a subdomain resulting from gluing together a collection of Giroux domains $G_{\Sigma_1},\dots, G_{\Sigma_k}$, we shall denote $\mathcal{O}(N)=\mathcal{O}(G_{\Sigma_1})\cap\dots\cap\mathcal{O}(G_{\Sigma_k})\subseteq H_2(M;\mathbb{R})$.

The following theorem implies Theorem \ref{AlgvsGiroux} (see Example \ref{GTDom} below): 

\begin{thm}\label{Girouxthen1torsion}
If a contact manifold $(M^{2n+1},\xi)$ admits a contact embedding of a subdomain $N$ obtained by gluing two Giroux domains $G_{\Sigma_-}$ and $G_{\Sigma_+}$, such that $\Sigma_+$ has a boundary component not touching $\Sigma_-$, then $(M^{2n+1},\xi)$ has $\Omega$-twisted algebraic 1-torsion, for every $\Omega \in \mathcal{O}(N)$. Moreover, it is also algebraically overtwisted if $N$
contains any blown down boundary components.
\end{thm}

The motivating example of an explicit model of such subdomain is the following:

\begin{example}\label{GTDom}
Consider $(Y,\alpha_+,\alpha_-)$ a Liouville pair on a closed manifold $Y^{2n-1}$. As in the introduction, consider the \emph{Giroux $2\pi$-torsion domain} modeled on $(Y,\alpha_+,\alpha_-)$, given by the contact manifold $(GT,\xi_{GT}):=(Y\times [0,2\pi]\times S^1,\ker\lambda_{GT})$, where 
\begin{equation}\label{lambdaGT}
\begin{split}
\lambda_{GT}&=\frac{1+\cos(r)}{2}\alpha_++\frac{1-\cos(r)}{2}\alpha_-+\sin (r)d\theta\\
\end{split}
\end{equation} and the coordinates are $(r,\theta)\in [0,2\pi] \times S^1$.

We may write $\lambda_{GT}=f(\alpha+d\theta),$ where $$f=f(r)=\sin (r),$$$$\alpha=\frac{1}{2}(e^{u(r)}\alpha_++e^{-u(r)}\alpha_-),$$$$u(r)=\log \frac{1+\cos(r)}{\sin (r)},$$ Here, $u:(0,\pi)\rightarrow \mathbb{R}$ is an orientation reversing diffeomorphism, whereas $u:(\pi,2\pi)\rightarrow \mathbb{R}$ is an orientation-preserving one, which is used to pull-back the Liouville form $\frac{1}{2}(e^u\alpha++e^{-u}\alpha_+)$ defined on $Y \times \mathbb{R}$ via a map of the form $v=id\times u$. 

This means that we may view $(GT,\xi_{GT})$ as being obtained by gluing two Giroux domains of the form $GT_-=Y \times [0,\pi] \times S^1$ and $GT_+=Y \times [\pi,2\pi] \times S^1$, along their common boundary, an $\xi_{GT}$-round hypersurface modeled on $(Y,\alpha_-)$. 

\vspace{0.5cm}

Therefore, from Theorem \ref{Girouxthen1torsion} we obtain Theorem \ref{AlgvsGiroux} as a corollary.

\end{example}

\subsection{Giroux SOBDs}
\label{GirouxSOBDs}

We consider a specially simple kind of spinal open book decompositions (SOBDs), which arise on manifolds which have been obtained by gluing a family of Giroux domains along a collection of common boundary components, each a round hypersurface modeled on some contact manifold. Such is the case of the Giroux $2\pi$-torsion domains $GT$. Basically, these SOBDs are obtained by declaring suitable collar neighbourhoods of each gluing hypersurface to be \emph{paper} components, whereas the \emph{spine} components are the complement of these neighbourhoods. They have the desirable features that the fibrations are trivial, and that they have 2-dimensional \emph{pages} (see \cite{Mo2} for definitions). 

\paragraph*{Construction of the SOBD.} Let $(G_\pm=\Sigma_\pm\times S^1, \xi_\pm=\ker(f_\pm(\beta_\pm+d\theta)))$ be two Giroux domains which one wishes to glue along a round-hypersurface $H=Y\times S^1$ modeled on some contact manifold $(Y,\xi_0)$ and lying in their common boundaries. Fix choices of collar neighbourhoods $\mathcal{N}_\pm(H)$ of $H$ inside $G_\pm$, of the form $Y \times S^1 \times [0,1]$. Take coordinates $s_\pm \in [0,1]$ so that $H=\{s_\pm=0\}$, and $\theta \in S^1$, such that the contact structures $\xi_\pm$ are given by the kernel of the contact forms $\beta_0+s_\pm d\theta$. Here, $\beta_0$ is a contact form for $\xi_0$. In these coordinates, $\beta_\pm=\beta_0/s_\pm$ and $f=f(s_\pm)=s_\pm$, and the corresponding \emph{ideal} Liouville vector fields are $V_\pm=-s_\pm \partial_{s_\pm}$. We glue $\mathcal{N}_\pm(H)$ together in the natural way, by taking a coordinate $s\in [-1,1]$, so that $s=-s_-$ and $s=s_+$, and $s=0$ corresponds to $H$. We thus obtain a collar neighbourhood $\mathcal{N}(H)=\mathcal{N}_+(H) \bigcup_\Phi \mathcal{N}_-(H)\simeq H\times [-1,1]$, where we denote by $\Phi$ the resulting gluing map. Doing this for each of the boundary components that we glue together, we obtain a decomposition for $$M:=G_+\bigcup_{\Phi} G_-=\left(\Sigma_+\bigcup_{\Phi} \Sigma_-\right)\times S^1,$$ given by $$M=M_P \cup M_\Sigma,$$ where $M_P$ (the \emph{paper}) is the disjoint union of the collar neighbourhoods of the form $\mathcal{N}(H)$ for each of the gluing round-hypersurfaces $H$, and $M_\Sigma$ (the \emph{spine}) is the closure of its complement in $M$. We have a natural fibration structure $\pi_P$ on $M_P$, which is the trivial fibration over the disjoint union of the hypersurfaces $H$, so that the pages (the fibers of $\pi_P$) are identified with the annuli $P:=S^1\times [-1,1]$. We also have an $S^1$-fibration $\pi_\Sigma$ on $M_\Sigma$, which is also trivial. It has as base the disjoint union of $\Sigma_+$ and $\Sigma_-$ minus the collar neighbourhoods, which we denote by $\Sigma$. Let us assign a \emph{sign} $\epsilon(G_\pm):=\pm 1$ to each of the Giroux subdomains $G_\pm \subseteq M$.

We can also blow-down boundary components in the Giroux subdomains, and in this case we may extend our SOBD by declaring the $\mathbb{D}\times Y$ we glued in the blow-down operation to be part of the paper $M_P$ (so that its pages are disks). We declare the blown-down components to be part of the paper.  
\vspace{0.5cm}

We shall refer to the SOBDs obtained by the procedure described above as \emph{Giroux} SOBDs, where are also allowing the number of \emph{Giroux subdomains} involved (which is the same as the number of spine components) to be arbitrary. 

\vspace{0.5cm}

We will also fix collar neighbourhoods of the components of $\partial M$, which correspond to the non-glued hypersurfaces, in the very same way as we did before for the glued ones. Their components look like $\mathcal{N}_\Sigma(H):=H\times [0,1] = Y\times S^1\times [0,1]$, for some non-glued hypersurface $H$. There is a coordinate $s\in [0,1]$, such that $\partial M=\{s=0\}$, and such that the contact structure on the corresponding Giroux domain is given by $\ker(\beta_0+s d\theta)$, for some contact form $\beta_0$ for $(Y,\xi_0)$. Denote by $\mathcal{N}_\Sigma$ the disjoint union, over all unglued $H$'s, of all of the $\mathcal{N}_\Sigma(H)$'s. 

\paragraph*{The Giroux form.} If we denote by $\lambda_\pm= f_\pm(\beta_\pm + d\theta)$, we can extend the expression $\lambda_+=\beta_0+sd\theta$, a priori valid on $\mathcal{N}_+(H)$, to $\mathcal{N}(H)$, by the same formula. Observe that, by choice of our coordinates, the resulting 1-form glues smoothly to $\overline{\lambda}_-:=f_-(\beta_--d\theta)$. Therefore, we may still think of the contact structure $\overline{\xi}_-=\ker \overline{\lambda}_-$ as a \emph{contactization} contact structure over the region $(\Sigma_-\times S^1)\backslash\mathcal{N}_-(H)$, with the caveat that we need to switch the orientation in the $S^1$-direction. 

Denote the resulting contact form by $\lambda:=\lambda_+\cup_\Phi \overline{\lambda}_-$. We may globally write it as $\lambda=f(\beta+d\theta).$ Here, $f$ is a function which is either strictly positive or strictly negative over the interior of the Giroux domains, and vanishes precisely at the (glued and unglued) boundary components. The $1$-form $\beta$ coincides with the Liouville forms $\pm \beta_\pm$ where these are defined, and is undefined along said boundary components. In the case of blown-down components, the contact form $\lambda$ also extends in a natural way.

\begin{figure}[t]
\centering
\includegraphics[width=0.45\linewidth]{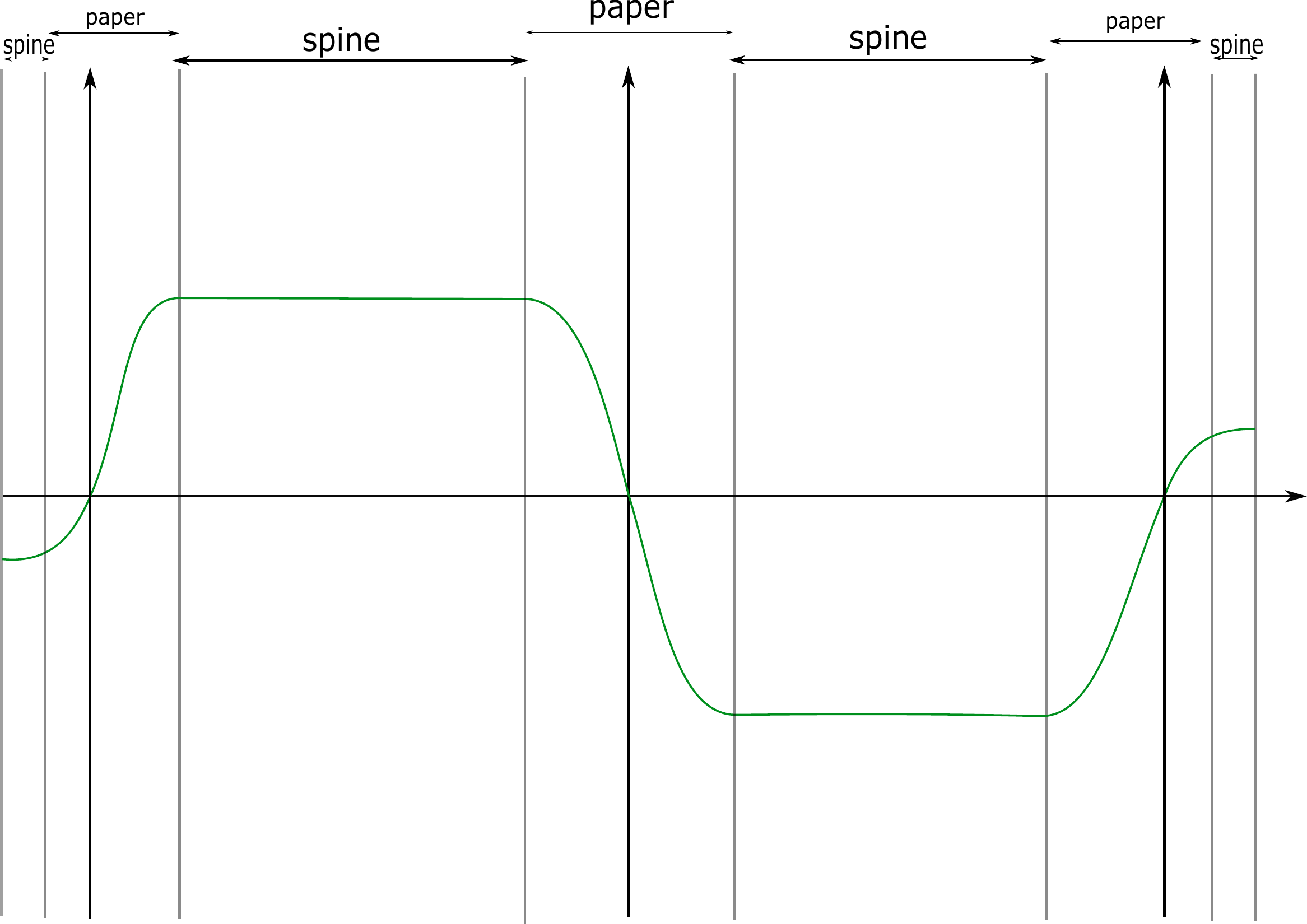}
\caption{\label{fGiroux} The isotoped function $f$ in the extended Giroux torsion domain $GT^\epsilon$.}
\end{figure}

From this construction, in $M$, each of the subdomains $G\subset M$ used in the gluing procedure carries a \emph{sign} $\epsilon(G)$, which we define as the sign of the function $f\vert_{\mathring{G}}$. Observe that $f$ coincides with $\epsilon(G)$ along every $H \times \{s_\pm=1\}$, the inner boundary components of the collar neighbourhoods. Therefore, we may isotope it relative every collar neighbourhood to a smooth function which is constant equal to $\epsilon(G)$ along $G\backslash \left((M^\delta_P \cap G) \cup \mathcal{N}^\delta_\partial\right)$ (see Figure \ref{fGiroux}). Here, $M^\delta_P$ and $\mathcal{N}^\delta_\partial$ denote small $\delta$-extensions in the interval direction of $M_P$ and $\mathcal{N}_\Sigma$, respectively, so that now $|s| \in [0,1+\delta]$. The regions $M_P^\delta\backslash M_P$ play the role of smoothened corners. Observe that isotoping $f$ as we just did does not change the isotopy class of $\lambda$, by Gray stability (version with boundary), and has the effect of transforming the Reeb vector field of $\lambda$ into $\epsilon(G)\partial_\theta$ along $M_Y$. Observe also that along the paper, we have $\lambda=\beta_0+\gamma$, where $\gamma=sd\theta$ is a Liouville form for $S^1\times [-1,1]$. Since this manifold is trivially a cylindrical Liouville semi-filling, we may view $\lambda$ as a Giroux form, as we did with contact form constructed in Section \ref{model}. 

\paragraph*{$\epsilon$-extension.} It will be convenient to consider an $\epsilon$-extension of our SOBD, which we call $M^\epsilon$, by gluing small collars to the boundary. To each component $\mathcal{N}_\Sigma(H)$ for $H$ lying in $\partial G$ for a subdomain $G$ of $M$, we glue a collar neighbourhood of the form $H \times [-\epsilon,0]$ for some small $\epsilon>0$. We extend our function $f$ so that the extended version coincides with $\epsilon(G)id$ near $s=0$, $f^\prime(-\epsilon)=0$, and $sgn(f^\prime(s))=\epsilon(G)\neq 0$ for $s>-\epsilon$ (see Figure \ref{fGiroux}). We will define the paper $M^\epsilon_P$ to be the union of $M_P$ and the region $H \times [-\epsilon/2,1+\delta]$, and the spine $M_\Sigma^\epsilon$ to be the union of $M_\Sigma$ with $H \times [-\epsilon,-\epsilon/2]$. The point of this extension is that now the Reeb vector field of the extended $\lambda$ coincides with $-\epsilon(G)\partial_\theta$ along the boundary. 

\vspace{0.5cm}

\begin{remark}
Because of the above discussion on orientations, from which we gathered that the $S^1$-orientation that we need depends on the sign of the Giroux domain, we rule out the case where we have a sequence $(G_0,\dots,G_{N-1})$ of Giroux subdomains of $M$, where $G_i$ has been glued to $G_{i+1}$ (modulo $N$) along some collection of boundary components, and $N>1$ is odd. This condition is to be taken as part of the definition of a Giroux SOBD.

\end{remark}

\subsection{Prequantization SOBDs (non-trivial case)} \label{PrequantSOBDs}

We can generalize the previous construction to the case where the $S^1$-bundles are not necessarily globally trivial, but trivial on the boundary components which are glued together.

\begin{defn}\cite{DiGe}\label{IdealLDomDIGE} Let $\Sigma$ be a compact $2n$-manifold with boundary, $\omega_{int}$ a symplectic form on $\mathring{\Sigma}$, $\xi_0$ a contact structure on $\partial \Sigma$, and $c$ a cohomology class in $H^2
(\Sigma; \mathbb{Z})$. The tuple $(\Sigma, \omega_{int}, \xi_0, c)$ is called an \emph{ideal Liouville
domain} if for some (and hence any) closed 2-form $\omega$ on $\Sigma$ with $-[\omega/2\pi] = c \otimes \mathbb{R} \in H^2(\Sigma;\mathbb{R})$ there exists a 1-form $\beta$ on $\mathring{\Sigma}$ such that 
\begin{itemize}
\item $d\beta = \omega_{int} - \omega$ on $\mathring{\Sigma}$.
\item For any smooth function $f:\Sigma \rightarrow [0,\infty)$ with regular level set $\partial\Sigma=f^{-1}(0)$, the 1-form $\beta_0=f\beta$ extends smoothly to $\partial \Sigma$ such that its restriction to $\partial \Sigma$ is a contact form for $\xi_0$.
\end{itemize} 
\end{defn}
 
An ideal Liouville domain $(\Sigma, \omega_{int}, \xi_0, 0)$ is an ideal Liouville domain in the sense of Giroux, where we may take $\omega = 0$. In the case where $c\vert_{\partial \Sigma}=0$, so that we may take $\omega=0$ in any collar neighbourhood  $[0,1)\times \partial \Sigma$ of the boundary, we will call $(\Sigma, \omega_{int}, \xi_0, c)$ an \emph{ideal strong symplectic filling}.

\vspace{0.5cm}

Now let $\pi: M\rightarrow \Sigma$ be the principal $S^1
$-bundle over $\Sigma$ of (integral) Euler class $c$.
Choose a connection 1-form $\psi$ with curvature form $\omega$ on this bundle. As in the trivial case $c=0$ (where we may take $\psi=d\theta$), the 1-form $f(\psi + \beta)=f\psi + \beta_0$ defines a contact structure $\xi$ on $M$, and we call $(M, \xi)$ the \emph{contactization} of $(\Sigma, \omega_{int}, \xi_0, c)$. Observe that $d(\psi + \beta)=\pi^*(\omega + \omega_{int}-\omega)=\pi^*\omega_{int}$, where we identify $\beta$ with $\pi^*\beta$, so that on $\mathring{\Sigma}$, $\xi$ may be regarded as the prequantization contact structure corresponding to $\omega_{int}$. 

\begin{defn}\cite{DiGe} Let $B$ be a closed, oriented manifold of dimension $2n$. An ideal Liouville splitting of class $c \in H^2(B; \mathbb{Z})$ is a decomposition $B = B_+ \bigcup_\Gamma B_-$ along a two-sided (but not necessarily connected) hypersurface $\Gamma$, oriented as the boundary of $B_+$, together with a contact structure $\xi_0$ on $\Gamma$ and symplectic forms $\omega_{\pm}$ on $\pm\mathring{B}_\pm$, such that $(\pm B_\pm, \omega_\pm, \xi_0, \pm c\vert_{B_\pm})$
are ideal Liouville domains.
\end{defn}

Proposition 3.6 in \cite{DiGe} tells us that an $S^1$-invariant
contact structure $\xi = \ker(\beta_0 + f\psi)$ on the principal $S^1$-bundle $\pi: M \rightarrow B$ defined by $c \in H^2(B; \mathbb{Z})$ leads to an ideal Liouville splitting of $B$ of class $c$, along the dividing set $$\Gamma=\{b \in B: \partial_\theta\vert_b \in \xi\}=f^{-1}(0)$$ 

If we require that $c \in H^2(B,\Gamma;\mathbb{Z})$, so that $c\vert_{\Gamma} =0$, then the $S^1$-principal bundle $\pi$ is trivial in any collar neighbourhoods $[-1,1]\times \Gamma$ of $\Gamma$, along which we may take $\psi=d\theta$, $\omega=0$. In this situation, we may regard the total space $M$ as carrying what we will call a \emph{prequantization SOBD}, which is obtained in the same way as we defined a Giroux SOBD. That is, we declare suitable collar neighbourhoods of $\Gamma$ to be paper components, which we can do by the triviality assumption on $c$ along $\Gamma$. The only difference now is that the spine is no longer globally a trivial $S^1$-bundle. Reciprocally, one can glue prequantization spaces over ideal strong symplectic fillings to obtain a manifold $M$ (possibly with non-empty boundary) carrying an $S^1$-invariant contact structure $\xi=\ker(\beta_0+f\psi)$, and a prequantization SOBD. The gluing construction is the same as for Giroux SOBDs, and is a particular case of the gluing construction of Thm.\ 4.1 in \cite{DiGe}.  

Observe that, using the coordinates of the previous section, the contact $1$-form $\alpha=\beta_0+f\psi$ defining $\xi$ satisfies that $d\alpha=ds \wedge d\theta>0$ is positive along the pages $S^1\times [-1,1]$, and, after isotoping $f$ in the same way as for Giroux SOBDs, its Reeb vector field is tangent to the $S^1$-fibers of $\pi_\Sigma$ along $M_\Sigma$. Moreover, along the boundary, it induces the distribution $\ker \alpha \cap T(\partial M)=\xi_0 \oplus TS^1$, with characteristic foliation given by the $S^1$-direciton, and its Reeb vector field coincides with $R_0$, the Reeb vector field of $\beta_0$. From this, and having the definition in \cite{LVHMW}, and the standard one due to Giroux, both in mind, one can define

\begin{defn}
A \emph{Giroux form} for a prequantization SOBD (in particular, for a Giroux SOBD) is any contact form inducing a contact structure which is isotopic to the $S^1$-invariant contact structure $\xi=\ker(\beta_0+f\psi)$.
\end{defn} 

With this definition, any $S^1$-invariant contact form $\alpha$ inducing an $S^1$-invariant contact structure in a principal bundle defined by a cohomology class $c$ satisfying $c\vert_{\Gamma}=0$, where $\Gamma$ is a set of dividing hypersurfaces for $\alpha$, is a Giroux form for the induced prequantization SOBD.

\subsection{Proof of Theorem \ref{Girouxthen1torsion}}

The proof of this theorem is a reinterpretation of what we did in the construction of our contact manifold models of Section \ref{model}. 

\begin{proof}[Thm. (\ref{Girouxthen1torsion})]

Let $N$ be a subdomain as in the hypothesis, carrying a Giroux form $\lambda=f(\beta+d\theta)$ obtained by gluing. Since the contact embedding condition is open, and we are assuming that there are boundary components of $\Sigma_+$ not touching $\Sigma_-$, we can find a small $\epsilon>0$ such that $M$ admits a contact embedding of the $\epsilon$-extension $N^\epsilon$. Endow this extension with a Giroux SOBD $N^\epsilon=M^\epsilon_P \cup M^\epsilon_\Sigma$ as in the previous section. On this decomposition, add corners where spine and paper glue together as we explained before. Also choose a small Morse/Morse--Bott function $H$ in $M_\Sigma^\epsilon$, which lies in the isotopy class of $f$, and vanishes as we get close to the paper. With this data, we then may construct a Morse/Morse--Bott contact form $\Lambda$ on $N$ which lies in the isotopy class of $\lambda$, along with a SHS deforming it, in the analogous way as done in Section \ref{model}. 

Since in our situation we are allowing $\Sigma_\pm$ to be more general than a semi-filling $Y \times I$, we need to specify what we mean by the Morse--Bott situation. We will take our Morse--Bott function $H$ so that it depends only on the interval parameter along the collar neighbourhoods $\mathcal{N}_\Sigma$ close to the boundary and along a slightly bigger copy of $M_P^\epsilon$, and matches the function $f$ close to $\partial N$. In particular, $\partial N$ is a Morse--Bott submanifold. We also impose that $H$ is Morse in the interior of the components of $M_\Sigma$ which are away from the boundary. For simplicity, we will assume that, besides the boundary, $H$ only has exactly one Morse--Bott submanifold close to each boundary component of $M_\Sigma$ which is glued to a paper component, of the form $Y\times \{t\}$ for some $t$ (see Figure \ref{Girouxd}). The Morse situation is then obtained by a perturbation of this situation obtained by choosing Morse functions along the Morse--Bott submanifolds. Observe that we may always choose the interior Morse-Bott submanifolds so that they lie in the cylindrical ends of the Liouville domains $\Sigma_\pm$.

\begin{figure}[t]
\centering
\includegraphics[width=0.60\linewidth]{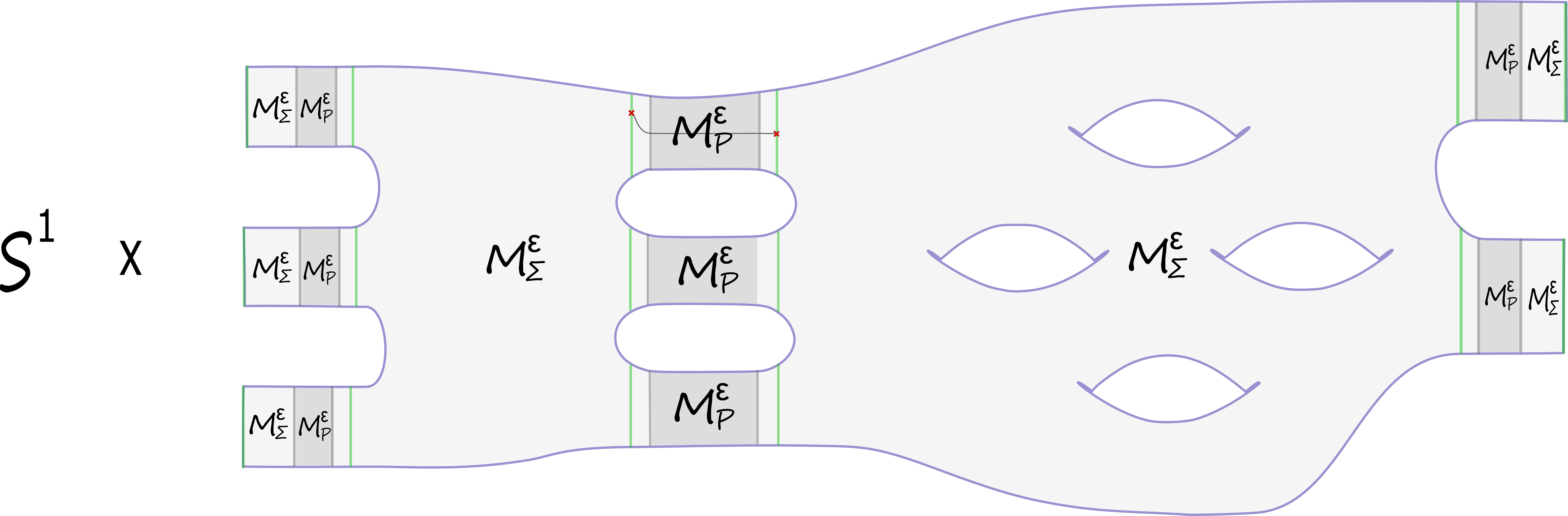}
\caption{\label{Girouxd} An $\epsilon$-extension of a domain $N$ consisting of two Giroux domains glued along three boundary components. We draw the Morse--Bott submanifolds of $H$ in green. These are a ``barrier'' for page-like holomorphic curves, needed to adapt Theorem \ref{uniqueness}. We also draw an index $1$ holomorphic cylinder (of the first kind).}
\end{figure} 

We have an almost complex structure $J^\prime$ compatible with the SHS, for which we get a stable finite energy foliation by holomorphic cylinders, which come in two types: either they are obtained by gluing constant lifts of the cylindrical pages and flow-line cylinders over $M^\epsilon_\Sigma$ along the corners (first kind); or they are flow-line cylinders completely contained in $M_\Sigma^\epsilon$ (second kind). Both kinds have as asymptotics simply covered Reeb orbits $\gamma_p$ corresponding to critical points $p$ of $H$, either along the boundary, or at interior points of $\Sigma$ (see Figure \ref{foliationGirouxTorsion}). These cylinders have two positive ends if its corresponding critical points lie in Giroux subdomains with different sign, and one positive and one negative end if these signs agree. The Fredholm index formula is exactly as before. 

One can adapt the proof of Fredholm regularity in the Morse-Bott case (Section \ref{Regularity}) to this particular situation, to obtain regularity for cylinders of both kinds. There is no change in the proof for cylinders of the second kind. For regularity for cylinders of the first kind, we adapt the proof of Section \ref{reggenuszero}. For this, we observe that our choices of interior Morse-Bott submanifolds imply that the cylindrical ends of the cylinders of the first kind lie completely over cylindrical ends of the Liouville domains $\Sigma_\pm$ (see Figure \ref{Girouxd}). This means that we may assume that the obvious analogous version of Assumptions \ref{Assumptions} hold, and the rest of the proof is a routine adaptation. Fredholm regularity in the sufficiently nearby Morse case follows from the implicit function theorem.   

\begin{figure}[t]
\centering
\includegraphics[width=0.50\linewidth]{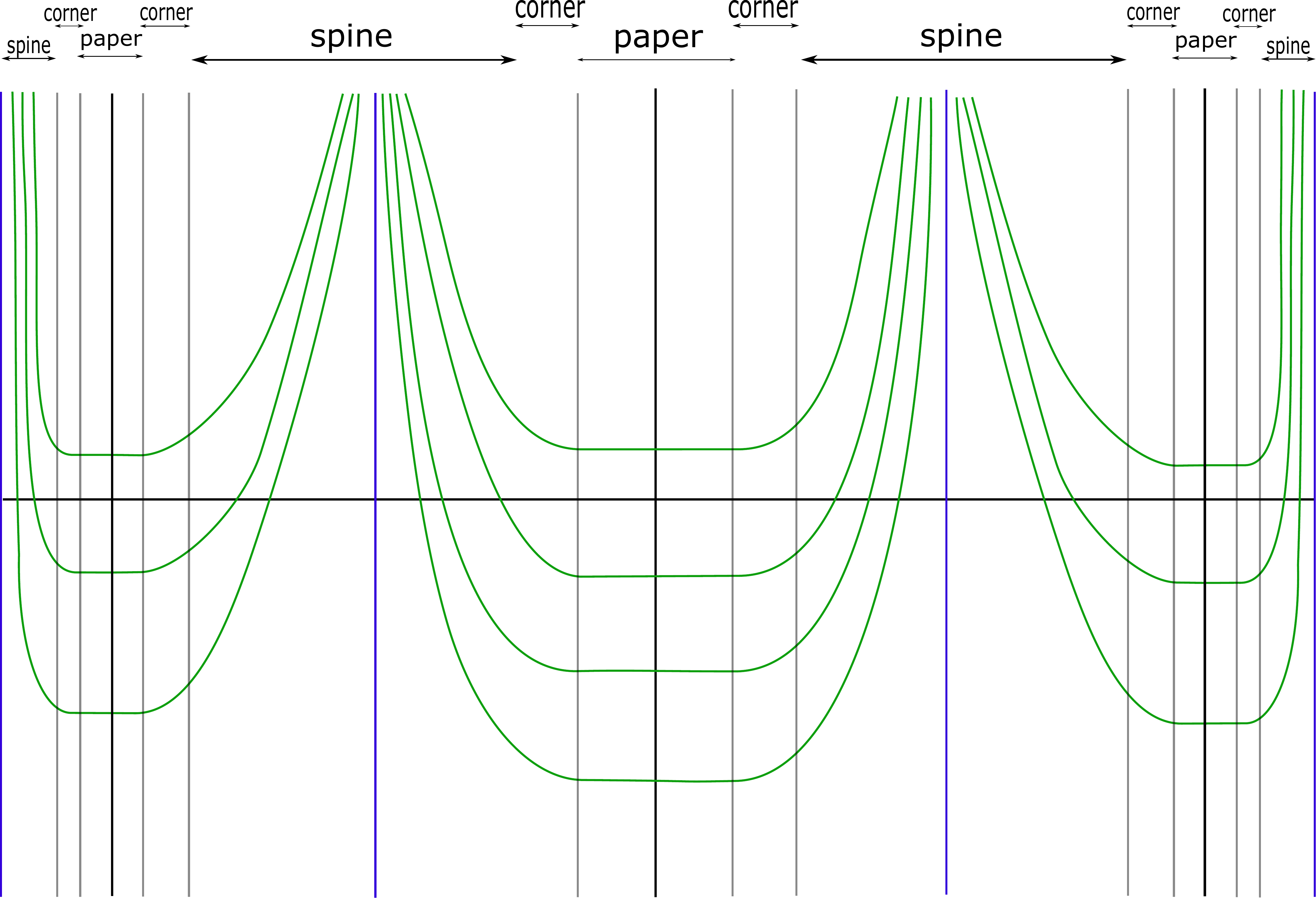}
\caption{\label{foliationGirouxTorsion} The foliation by holomorphic cylinders (with respect to SHS data) of the symplectization of $GT^\epsilon$.}
\end{figure}

We also have a version of the uniqueness Theorem \ref{uniqueness}, adapted to this situation. We have a few small changes, as follows. In the hypothesis we require that the positive asymptotics are simply covered; correspond to critical points all of which, in the Morse--Bott case, belong to a Morse--Bott submanifold of the type already described; and any two lie in distinct components of $M_\Sigma$, both of which are away from the boundary, separated by a paper component. The number $k$ is replaced by $2$ in the proof. We have more than just one Morse--Bott submanifold now, but, for case B where $u$ lies over some paper component $P=Y \times I \times S^1$, we can still get a similar upper bound on the energy:

\begin{equation}
\begin{split}
\mathbf{E}(u)\leq & 2\pi\left(\sum_{z \in \Gamma^+}e^{\epsilon H(p_z)}- \sum_{z \in \Gamma^-}T_z\right)\\
\leq &2\pi \sum_{z \in \Gamma^+}e^{\epsilon H(p_z)} \\
\leq &4\pi ||e^{\epsilon H}||_{C^0},
\end{split}
\end{equation}
where $T_z>0$ is the action of the Reeb orbit corresponding to $z \in \Gamma^-$, which a priori might even lie in $M\backslash N$ (but not a posteriori). Here, the energy is computed with respect to a $2$-form which is the derivative of a contact form for $\xi$ outside of $N$, and which restricts to the $2$-form of the fixed SHS on $N$. Here we use that we have chosen the SHS so that it is contact near $\partial N$. 

Using that the number of boundary components of each paper component is 2, as before we get
$$
\mathbf{E}(u)\geq 4 \pi e^{-\delta}\mbox{deg}_P(u),
$$
for some small $\delta>0$. Here, we denote by $\mbox{deg}_P(u)$ the covering degree of the restriction to $u$ of the projection to $I \times S^1$, over the component $P$. The rest of the proof is the same. Observe that since we are assuming that the critical points appearing in the positive asymptotics are away from the boundary of $N$, and the only way of venturing into $M\backslash N$ is to escape through the latter, any holomorphic curve with these positive asymptotics which leaves $N$ would necessarily need to go through a paper component, so is dealt with by case B in the proof of \ref{uniqueness} (which deliberately does not need to assume that the whole curve stays over $N$).

\vspace{0.5cm}

Make all the necessary choices to have an SFT differential 
$$
\mathbf{D}_{SFT}:\mathcal{A}(\Lambda)[[\hbar]]\rightarrow \mathcal{A}(\Lambda)[[\hbar]],
$$
which computes $H_*^{SFT}(N,\ker \Lambda;\Omega\vert_N)$, where $\Omega \in \mathcal{O}(N)$. Extend these choices to $M$ so as to be able to compute $H_*^{SFT}(M,\xi;\Omega)$.

Let $e$ be a maximum (index $2n$) of $H$ in $M_\Sigma \cap (\Sigma_-\times S^1)$, and let $h$ be an index $1$ critical point of $H$ in $M_\Sigma \cap (\Sigma_+\times S^1)$. We take both to lie in the Morse--Bott manifolds of the Morse--Bott case, before the Morse perturbation. Denote by $q_e$, $q_h$ the corresponding SFT generators. Define 
$$
Q=q_hq_e,
$$
an element of $\mathcal{A}(\Lambda)$. There is a unique (perturbed) cylinder $u$ of the first kind which has $\gamma_e$ and $\gamma_h$ as positive asymptotics, with $\mbox{ind}(u)=1$. If we choose $H$ so that it does not have any minimum, by uniqueness any other holomorphic curve over $M$ which may contribute to the differential of $Q$ is a flow-line cylinder completely contained in $M_\Sigma$, connecting an index $2n-1$ critical point $p$ with $e$. 

Denote by $\overline{d}$ the element in $H_2(M;\mathbb{R})/\ker \Omega$ defined by any $d \in H_2(M;\mathbb{R})$, by $\mathcal{M}(H;p,e)$ the space of positive flow-lines connecting $p$ with $e$, $u_\gamma$ the flow-line cylinder corresponding to a flow-line $\gamma$, and $\epsilon(u)$ the sign of the holomorphic curve $u$ assigned by a choice of coherent orientations. Since $H_{2n}(\Sigma_\pm)=0$, elements of $\mathcal{M}(H;p,e)$ come in evil twins pairs $\gamma \leftrightarrow \overline{\gamma}$ such that $\epsilon(u_\gamma)=-\epsilon(u_{\overline{\gamma}})$. As in Section \ref{Torsion}, one can choose suitable spanning surfaces such that $\overline{[u_\gamma]}=\overline{[u_{\overline{\gamma}}]}$. 

Then,  
\begin{equation}
\begin{split}
\mathbf{D}_{SFT}Q=&z^{\overline{[u]}}\hbar+\sum_{\substack{\gamma \in \mathcal{M}(H;p,e)\\ \mbox{ind}_p(H)=2n-1}} \epsilon(u_\gamma)z^{\overline{[u_\gamma]}}q_pq_h\\
=& z^{\overline{[u]}}\hbar+\frac{1}{2}\sum_{\substack{\gamma \in \mathcal{M}(H;p,e)\\ \mbox{ind}_p(H)=2n-1}} \left(\epsilon(u_\gamma)+\epsilon(u_{\overline{\gamma}})\right) z^{\overline{[u_\gamma]}}q_pq_h\\
=&z^{\overline{[u]}}\hbar,\\
\end{split}
\end{equation}  
which proves that $(N,\ker \Lambda)$ has $\Omega\vert_N$-twisted 1-torsion. This implies that $(M,\xi)$ has $\Omega$-twisted 1-torsion, since our uniqueness theorem gives that there are no holomorphic curves with asymptotics in $N$ which venture into $M\backslash N$, and therefore $H_*^{SFT}(N,\ker \Lambda; \Omega\vert_N)$ embeds into $H_*^{SFT}(M,\xi; \Omega)$. Here, we use that $\ker \Lambda$ is isotopic to $\ker \lambda=\xi\vert_N$.

\vspace{0.5cm}

For the second statement, assume that we have a blown-down boundary component, so that the corresponding $\mathbb{D}\times Y$ is a paper component with disk pages. As in Section \ref{Curves}, we have that the disk pages lift as finite-energy holomorphic planes with a single positive asymptotic $\gamma_p$, corresponding to a critical point $p$ in $M_\Sigma$. If $u$ is such a plane, its index is $\mbox{ind}(u)=ind_p(H)$. Take $p$ so that $ind_p(H)=1$, and let $P=q_{\gamma_p}$. Since there is no minima for $H$, by our uniqueness theorem we know that
$$
\mathbf{D}_{SFT}\left(z^{-\overline{[u]}}P\right)=1
$$  
Since the choice of coefficients is arbitrary, $(N,\ker \Lambda)$ is algebraically overtwisted, which implies that $(M,\xi)$ is.

\end{proof}

In view of the definition of a prequantization SOBD, and the fact that we may still construct a foliation by flow-line holomorphic cylinders over a non-trivial prequantization space \cite{Mo2,Sie}, one has the following:

\begin{thm}
In dimension $3$, consider a principal $S^1$-bundle obtained by gluing a collection of prequantization spaces over ideal strong symplectic fillings, such that there are two of them, one with a boundary component not touching the other after the gluing. Then the resulting $S^1$-invariant contact structure has (untwisted) algebraic 1-torsion. 
\end{thm}

In dimension at least $5$, the only technical ingredient for the above result that is missing is the regularity of the holomorphic flow-line cylinders for sufficiently small $H$.

\newpage

\bibliographystyle{alpha}

\begin{thebibliography}{9}

\bibitem[BEM]{BEM} M. Borman, Y. Eliashberg, E.Murphy, Existence and classification of overtwisted contact structures in all dimensions. Acta Math. 215 (2015), no. 2, 281--361.

\bibitem[Bo02]{Bo}
F. Bourgeois. A Morse--Bott approach to contact homology. Symplectic and contact topology: interactions and perspectives (Toronto, ON/Montreal, QC, 2001), 55--77, Fields Inst. Commun., 35, Amer. Math. Soc., Providence, RI, 2003.

\bibitem[BoMon04]{BoMon}
F. Bourgeois, K. Mohnke, Coherent orientations in symplectic field theory. Math. Z. 248 (2004), no. 1, 123--146.


\bibitem[CL09]{CL09} K. Cieliebak, J. Latschev. The role of string topology in symplectic field theory. New perspectives and challenges in symplectic field theory, 113--146, CRM Proc. Lecture Notes, 49, Amer. Math. Soc., Providence, RI, 2009.

\bibitem[DiGe12]{DiGe} F. Ding, H. Geiges, Contact structures on principal circle bundles. Bull. Lond. Math. Soc. 44 (2012), no. 6, 1189--1202.

\bibitem[EGH00]{EGH} Y. Eliashberg, A. Givental, and H. Hofer, Introduction to symplectic field theory, Geom. Funct. Anal., Special Volume (2000), 560–673.

\bibitem[HT1]{HT1} M. Hutchings, C. Taubes.Gluing pseudoholomorphic curves along branched covered cylinders. I. J. Symplectic Geom. 5 (2007), no. 1, 43--137. 

\bibitem[HT2]{HT2} M. Hutchings, C. Taubes. Gluing pseudoholomorphic curves along branched covered cylinders. II. J. Symplectic Geom. 7 (2009), no. 1, 29--133.

\bibitem[vK]{vK} O. van Koert, Simple computations in prequantization bundles. availabe at http://www.math.snu.ac.kr/$\sim$ okoert/tools/CZ\_index\_BW\_bundle.pdf

\bibitem[KM-VHM-W]{KMvhMW} C. Kutluhan, G. Matic, J. Van Horn-Morris, and A. Wand, Algebraic torsion via Heegaard Floer homology. arXiv:1503.01685v3.

\bibitem[LW11]{LW} J. Latschev, C. Wendl, Algebraic Torsion in Contact manifolds, Geom. Funct. Anal. 21 (2011), no. 5, 1144-1195. 

\bibitem[L-VHM-W]{LVHMW} S. Lisi, J. Van Horn-Morris, C. Wendl. On symplectic fillings of spinal open book decompositions I: Geometric constructions. Preprint, arXiv:1810.12017.
 
\bibitem[Lutz77]{Lutz} R. Lutz, Structures de contact sur les fibr\`{e}s principaux en cercles de dimension trois,
Ann. Inst. Fourier (Grenoble) 27 (1977), no. 3, ix, 1–15 (French, with English summary).

\bibitem[MNW13]{MNW}
P. Massot, K. Niederkrüger, C. Wendl, Weak and strong fillability of higher dimensional contact manifolds, Inv. math., May 2013, 192 (2), 287-373.

\bibitem[McD91]{McD} D. McDuff, Symplectic manifolds with contact type boundaries, Inv. math., December 1991, 103 (1), 651-671.

\bibitem[Mit95]{Mit} Y. Mitsumatsu, Anosov flows and non-Stein symplectic manifolds, Annales de l’institut Fourier, tome 45, no 5 (1995), p. 1407-1421.

\bibitem[Mo]{Mo} A. Moreno, SFT computations and intersection theory in higher-dimensional contact manifolds. To appear.

\bibitem[Mo2]{Mo2} A. Moreno, Algebraic torsion in higher-dimensional contact manifolds (PhD thesis). Available at arXiv:1711.01562. 

\bibitem[MS19]{MS} A. Moreno, R. Siefring. \emph{Holomorphic curves in the presence of holomorphic hypersurface foliations}. Preprint. ArXiv:1902.02700.

\bibitem[Mori09]{Mori09} A. Mori, Reeb foliations on $S^5$ and contact $5$–manifolds violating the Thurston-Bennequin inequality, preprint, 2009, arXiv:0906.3237.

\bibitem[SZ92]{SZ} D.Salamon, E. Zehnder. Morse theory for periodic solutions of Hamiltonian systems and the Maslov index. Comm. Pure Appl. Math. 45 (1992), no. 10, 1303--1360. 

\bibitem[Sie11]{Sie11} R. Siefring. Intersection theory of punctured pseudoholomorphic curves. Geom. Topol. 15 (2011), no. 4, 2351--2457. 

\bibitem[Sie16]{Sie} R. Siefring, Finite-energy pseudoholomorphic planes with multiple asymptotic limits. Math. Ann. 368 (2017), no. 1-2, 367--390. MR3651577.

\bibitem[Sie]{Sie2} R. Siefring. Symplectic field theory and stable hamiltonian submanifolds: Intersection theory. In preparation.

\bibitem[Wen1]{Wen1} 
C. Wendl, Automatic transversality and orbifolds of punctured holomorphic curves in dimension
four, Comment. Math. Helv. 85 (2010), no. 2, 347–407.

\bibitem[Wen2]{Wen2}
C. Wendl, A hierarchy of local symplectic filling obstructions for contact 3-manifolds. Duke Math. J. 162 (2013), no. 12, 2197--2283.

\bibitem[Wen3]{Wen3} C. Wendl, Lectures on Symplectic Field Theory. arXiv:1612.01009v2.

\bibitem[Wen4]{Wen5} C. Wendl, Open book decompositions and stable Hamiltonian structures. Expo. Math. 28 (2010), no. 2, 187--199.

\bibitem[Wen5]{Wen8} C. Wendl, Contact 3-manifolds, holomorphic curves and intersection theory. Lecture notes. 	arXiv:1706.05540










          
\end{thebibliography}

\end{document}